\documentclass[11pt]{amsart}

\usepackage {amsmath, amssymb, amscd, mathrsfs,  amsthm, stmaryrd,  bbm, diagbox, enumerate, slashed, graphicx, color, subfig, comment, enumitem}
\usepackage[all, cmtip]{xy}
\usepackage[text={6.1in,8.6in},centering,letterpaper,dvips]{geometry}
\usepackage{url, hyperref}
\usepackage{amscd}
\usepackage{mathtools} 
\usepackage{amstext} 
\usepackage{array}  
\usepackage{xspace}

\newtheorem{thm}{Theorem}[section]
\newtheorem{theorem}[thm]{Theorem}

\newtheorem{corollary}[thm]{Corollary}

\newtheorem{conjecture}[thm]{Conjecture}

\newtheorem{question}[thm]{Question}
\theoremstyle{definition}

\newtheorem{definition}[thm]{Definition}
\theoremstyle{remark}

\newtheorem{remark}[thm]{Remark}

\newcommand\eps{\varepsilon}

\newcommand\Z{\mathbb{Z}}
\newcommand\R{\mathbb{R}}
\newcommand\Q{\mathbb{Q}}

\newcommand\Ta{\mathbb{T}_\alpha}
\newcommand\Tb{\mathbb{T}_\beta}

\def\Sym{\mathrm{Sym}}

\def\Z {{\mathbb{Z}}}
\def\R {{\mathbb{R}}}
\def\C {{\mathbb{C}}}
\def\Q {{\mathbb{Q}}}

\def\Link {\mathbb{L}}

\def\kk {{\mathbf{k}}}
\def\del {{\partial}}

\def\C {{\mathbb{C}}}

\def\spc {{\operatorname{Spin^c}}}
\newcommand\SpinC\spc

\def\rk {{\operatorname{rk} \hskip2pt}}
\def\fin\qedhere
\def\pr {{\text{pr}}}

\def\Kh{\mathit{Kh}}

\def\x{\mathbf x}
\def\y{\mathbf y}

\def\K {\mathcal K}

\def\N {\mathbb{N}}

\def\CP {\mathbb{CP}}
\def\bCP {\overline{\CP^2}}

\def\fin\qedhere
\def\kk{\Bbbk}

\def\Sz{\mathcal{S}}
\newcommand\HF{\mathit{HF}}
\newcommand\HFhat{\widehat{\HF}}
\newcommand\HFm{\HF^-}
\newcommand\HFp{\HF^+}
\newcommand\HFinfty{\HF^\infty}

\newcommand\HFKhat{\widehat{\mathit{HFK}}}
\newcommand\HFKt{\widetilde{\mathit{HFK}}}
\def\KBSM{\mathit{KBSM}}
\def\K{\mathbb{K}}

\title[From knots to four-manifolds]{From knots to four-manifolds}

\author[Ciprian Manolescu]{Ciprian Manolescu}
\address{Department of Mathematics, Stanford University\\
450 Jane Stanford Way, Building 380, Stanford, CA 94305-2125, USA.}
\email{cm5@stanford.edu}

\begin{document}

\begin{abstract} This is a survey article about the connections between knot theory and four-dimensional topology. Every four-manifold can be represented in terms of a link, by a Kirby diagram. This point of view has led to progress in computing invariants of smooth four-manifolds that can detect exotic structures. We explain  how this was done in two contexts: Heegaard Floer theory and skein lasagna modules. We also describe a program to understand four-manifolds through the properties of knots on their boundaries. 
\end{abstract}
\maketitle

\section{Introduction}
Knots and four-manifolds are two central notions in low-dimensional topology. The goal of this article is to explore some recent advancements that came out of the interplay between these notions.

Knot theory (briefly reviewed in Section~\ref{sec:knots}) is the study of closed loops in three-dimensional space. Knots can be distinguished by various algebraic invariants. Traditionally, such invariants took the form of numbers or polynomials; two famous ones are the Alexander and the Jones polynomial \cite{Alexander, Jones}. Since the turn of the millennium, an important development has been the emergence of homological invariants, which are multi-graded abelian groups. Examples include  knot Floer homology \cite{Knots, RasmussenThesis} and Khovanov homology \cite{Khovanov}. 

Modern four-dimensional topology has its origins in the early 1980's, when the independent work of Michael Freedman  and Simon Donaldson  revolutionized the field and made clear the distinction between topological and smooth four-manifolds \cite{Freedman, Donaldson}. Further progress came in the 1990's with the introduction of the Seiberg-Witten invariants \cite{Witten}. We will review a few highlights of four-dimensional topology in Section~\ref{sec:4m}, focusing on what are called {\em exotic pairs}: manifolds that are homeomorphic but not diffeomorphic.

Many of the early results on exotic pairs relied on solving PDE's using specific properties of the manifolds, such as having a K\"ahler form, a metric with positive scalar curvature, or a connected sum decomposition. The class of four-manifolds that admit these properties is limited, so it is desirable to have more topological methods that apply to {\em all} four-manifolds. In practice, topologists represent four-manifolds pictorially through Kirby diagrams. A Kirby diagram encodes a handle decomposition of the manifold, which is determined by a framed link in the connected sum of several copies of $S^1 \times S^2$; see  Section~\ref{sec:kirby}. Kirby diagrams are the key to relating four-manifold topology to the world of knots and links.

Can one compute four-manifold invariants by starting with a Kirby diagram as input? A success story in this direction is Heegaard Floer theory, which was developed by Ozsv\'ath and Szab\'o as a symplectic counterpart to Seiberg-Witten theory \cite{HolDisk, HolDiskTwo, HolDiskFour}. The resulting four-manifold invariants are conjecturally the same as the Seiberg-Witten invariants, and are algorithmically computable (at least mod $2$) from a Kirby diagram \cite{MOT}. While the general algorithm is not effective in practice, simpler methods can be used to compute specific examples, and they led to the discovery of new exotic pairs \cite{LLP}. We refer to Section~\ref{sec:HF} for more details.

Apart from computing existing four-manifold invariants, knot theory is also helpful in defining new such invariants. Skein lasagna modules \cite{MWW} are invariants associated to a framed link in the boundary of a four-manifold. Their definition relies on having a homological link invariant with good cobordism properties; Khovanov homology is an example. In particular, one can consider the skein lasagna module of the empty link in the boundary of a four-manifold, and hence obtain a four-manifold invariant. When the four-manifold is represented by a Kirby diagram, ways of computing the skein lasagna modules have been developed in \cite{MN, ManWW}. One remarkable application is the first analysis-free proof of the existence of exotic pairs of compact, orientable four-manifolds (with boundary)  \cite{RenWillis}. This circle of ideas is explored in  Section~\ref{sec:lasagna}.
 
Finally, in Section~\ref{sec:probing} we describe a more speculative program to construct exotic pairs by understanding which knots on the boundary of a four-manifold bound embedded disks in that manifold. There is some hope that this method could be used to tackle the smooth four-dimensional Poincar\'e conjecture, or at least detect exotic pairs of simply-connected four-manifolds with definite intersection form \cite{FGMW, MP}.

\medskip
 \textbf{Acknowledgments.} I am very grateful to Qianhe Qin for help in drawing many of the figures for this article, and Fabian Ruehle for Figure~\ref{fig:slice}. I would also like to thank Qiuyu Ren, Fabian Ruehle and Paul Wedrich  for comments that improved the paper, and Nathan Dunfield and Sherry Gong for discussions about their work \cite{DunfieldGong}.  I was partially supported by a Simons Investigator Award and the Simons Collaboration
Grant on New Structures in Low-Dimensional Topology.

\section{Knot theory}
\label{sec:knots}
A {\em link} is a collection of closed strings in three-dimensional space. More formally, it is a compact, smooth one-dimensional submanifold of $\R^3$. It is customary to equip the link with an orientation. Further, since topologists like compact spaces, it is common to add a point at infinity to $\R^3$, and view the link as an oriented submanifold of the three-dimensional sphere $S^3 = \R^3 \cup \{\infty\}.$

\begin{figure} \centering
{
\fontsize{9pt}{11pt}\selectfont
   \def\svgwidth{4.9in} 
  %% Creator: Inkscape 1.3.2 (091e20e, 2023-11-25), www.inkscape.org
%% Accompanies image file 'knots.pdf' (pdf, eps, ps)
\begingroup%
  \makeatletter%
  \providecommand\color[2][]{%
    \errmessage{(Inkscape) Color is used for the text in Inkscape, but the package 'color.sty' is not loaded}%
    \renewcommand\color[2][]{}%
  }%
  \providecommand\transparent[1]{%
    \errmessage{(Inkscape) Transparency is used (non-zero) for the text in Inkscape, but the package 'transparent.sty' is not loaded}%
    \renewcommand\transparent[1]{}%
  }%
  \providecommand\rotatebox[2]{#2}%
  \newcommand*\fsize{\dimexpr\f@size pt\relax}%
  \newcommand*\lineheight[1]{\fontsize{\fsize}{#1\fsize}\selectfont}%
  \ifx\svgwidth\undefined%
    \setlength{\unitlength}{553.48935105bp}%
    \ifx\svgscale\undefined%
      \relax%
    \else%
      \setlength{\unitlength}{\unitlength * \real{\svgscale}}%
    \fi%
  \else%
    \setlength{\unitlength}{\svgwidth}%
  \fi%
  \global\let\svgwidth\undefined%
  \global\let\svgscale\undefined%
  \makeatother%
  \begin{picture}(1,0.23408074)%
    \lineheight{1}%
    \setlength\tabcolsep{0pt}%
    \put(0,0){\includegraphics[width=\unitlength,page=1]{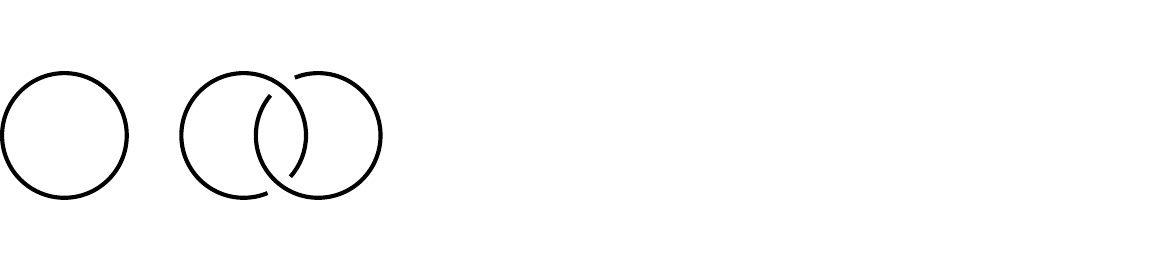}}%
    \put(0.03262551,0.00346889){\color[rgb]{0,0,0}\makebox(0,0)[lt]{\lineheight{1.25}\smash{\begin{tabular}[t]{l}$(a)$\end{tabular}}}}%
    \put(0,0){\includegraphics[width=\unitlength,page=2]{knots.pdf}}%
    \put(0.22040284,0.00346889){\color[rgb]{0,0,0}\makebox(0,0)[lt]{\lineheight{1.25}\smash{\begin{tabular}[t]{l}$(b)$\end{tabular}}}}%
    \put(0.42315576,0.00346889){\color[rgb]{0,0,0}\makebox(0,0)[lt]{\lineheight{1.25}\smash{\begin{tabular}[t]{l}$(c)$\end{tabular}}}}%
    \put(0.62590761,0.00346889){\color[rgb]{0,0,0}\makebox(0,0)[lt]{\lineheight{1.25}\smash{\begin{tabular}[t]{l}$(d)$\end{tabular}}}}%
    \put(0.86729827,0.00346889){\color[rgb]{0,0,0}\makebox(0,0)[lt]{\lineheight{1.25}\smash{\begin{tabular}[t]{l}$(e)$\end{tabular}}}}%
  \end{picture}%
\endgroup%

}
\caption{(a) the unknot; (b) the Hopf link; (c) the trefoil; (d) the Borromean rings; (e) the Conway knot}
\label{fig:knots}
\end{figure}

A connected link is called a {\em knot}. Some examples of knots and multi-component links are shown in Figure~\ref{fig:knots}. The main goal of knot theory is to classify knots (and links) up to isotopy, i.e., up to smooth deformation inside the ambient space through smooth submanifolds. The reader can easily find tables of knots and links online \cite{KnotAtlas, KnotInfo}.

\begin{figure} \centering
{
\fontsize{9pt}{11pt}\selectfont
   \def\svgwidth{5.1in} 
   %% Creator: Inkscape 1.3.2 (091e20e, 2023-11-25), www.inkscape.org
%% Accompanies image file 'reidemeister.pdf' (pdf, eps, ps)
\begingroup%
  \makeatletter%
  \providecommand\color[2][]{%
    \errmessage{(Inkscape) Color is used for the text in Inkscape, but the package 'color.sty' is not loaded}%
    \renewcommand\color[2][]{}%
  }%
  \providecommand\transparent[1]{%
    \errmessage{(Inkscape) Transparency is used (non-zero) for the text in Inkscape, but the package 'transparent.sty' is not loaded}%
    \renewcommand\transparent[1]{}%
  }%
  \providecommand\rotatebox[2]{#2}%
  \newcommand*\fsize{\dimexpr\f@size pt\relax}%
  \newcommand*\lineheight[1]{\fontsize{\fsize}{#1\fsize}\selectfont}%
  \ifx\svgwidth\undefined%
    \setlength{\unitlength}{305.82103098bp}%
    \ifx\svgscale\undefined%
      \relax%
    \else%
      \setlength{\unitlength}{\unitlength * \real{\svgscale}}%
    \fi%
  \else%
    \setlength{\unitlength}{\svgwidth}%
  \fi%
  \global\let\svgwidth\undefined%
  \global\let\svgscale\undefined%
  \makeatother%
  \begin{picture}(1,0.10437728)%
    \lineheight{1}%
    \setlength\tabcolsep{0pt}%
    \put(0.0643339,0.07411496){\makebox(0,0)[lt]{\lineheight{1.25}\smash{\begin{tabular}[t]{l}R1\end{tabular}}}}%
    \put(0,0){\includegraphics[width=\unitlength,page=1]{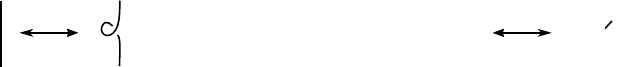}}%
    \put(0.79600244,0.07391469){\makebox(0,0)[lt]{\lineheight{1.25}\smash{\begin{tabular}[t]{l}R3\end{tabular}}}}%
    \put(0,0){\includegraphics[width=\unitlength,page=2]{reidemeister.pdf}}%
    \put(0.40224543,0.07377996){\makebox(0,0)[lt]{\lineheight{1.25}\smash{\begin{tabular}[t]{l}R2\end{tabular}}}}%
    \put(0,0){\includegraphics[width=\unitlength,page=3]{reidemeister.pdf}}%
  \end{picture}%
\endgroup%

}
\caption{Reidemeister moves}
\label{fig:reidemeister}
\end{figure}

Instead of considering smooth links, one can equivalently work with piecewise linear links, made of finitely many  non-intersecting segments connected at their ends. One can further arrange that the vertices of the link are on the three-dimensional integer lattice. It is then clear that the link is specified by a finite amount of data. This gives knot theory a strong combinatorial flavor. In practice, links are usually represented by their {\em planar diagrams} (regular projections to two-dimensional space), and then these projections are encoded in various ways (Gauss codes, DT codes, PD codes, etc.) Two planar diagrams represent the same link (up to isotopy) if and only if they can be related by a sequence of the {\em Reidemeister moves} shown in Figure~\ref{fig:reidemeister}.

To show that two links are isotopic, it suffices to exhibit a sequence of moves relating their diagrams. To show that two links (for example, the unknot and the trefoil) are not isotopic is more difficult: one needs an invariant, that is some quantity or algebraic object, typically defined in terms of a planar diagram, that is unchanged by Reidemeister moves. If the values of the invariant are different, we conclude that the two links are not isotopic. 

One of the best-known link invariants is the {\em Alexander polynomial} \cite{Alexander}. Given a link $L$, its Alexander polynomial $\Delta(L)$ is a Laurent polynomial in a variable $q^{1/2}$. (In the case of a knot, we only see integer powers of $q$.) There are several equivalent definitions of this polynomial; some of them can be found in the textbook \cite{LickorishBook}. The topological origin of the polynomial has to do with the homology of the infinite cyclic cover of the link complement. For our purposes, it suffices to say that the invariant is characterized by the normalization $\Delta(U)=1$ (where $U$ is the unknot) together with the skein relation
$$ \Delta(L_+) - \Delta(L_-) = (q^{1/2} - q^{-1/2}) \Delta(L_o),$$
which relates the polynomials of three links that differ at a single crossing in the manner shown in Figure~\ref{fig:orientedskein}. These properties allow for the Alexander polynomial to be calculated recursively, by starting from any planar diagram. For example, for the unlink of two components we have
$$ (q^{1/2} - q^{-1/2})\Delta( \scalebox{0.15}{%% Creator: Inkscape 1.4.2 (ebf0e940, 2025-05-08), www.inkscape.org
%% Accompanies image file 'unlink.pdf' (pdf, eps, ps)
\begingroup%
  \makeatletter%
  \providecommand\color[2][]{%
    \errmessage{(Inkscape) Color is used for the text in Inkscape, but the package 'color.sty' is not loaded}%
    \renewcommand\color[2][]{}%
  }%
  \providecommand\transparent[1]{%
    \errmessage{(Inkscape) Transparency is used (non-zero) for the text in Inkscape, but the package 'transparent.sty' is not loaded}%
    \renewcommand\transparent[1]{}%
  }%
  \providecommand\rotatebox[2]{#2}%
  \newcommand*\fsize{\dimexpr\f@size pt\relax}%
  \newcommand*\lineheight[1]{\fontsize{\fsize}{#1\fsize}\selectfont}%
  \ifx\svgwidth\undefined%
    \setlength{\unitlength}{118.5636205bp}%
    \ifx\svgscale\undefined%
      \relax%
    \else%
      \setlength{\unitlength}{\unitlength * \real{\svgscale}}%
    \fi%
  \else%
    \setlength{\unitlength}{\svgwidth}%
  \fi%
  \global\let\svgwidth\undefined%
  \global\let\svgscale\undefined%
  \makeatother%
  \begin{picture}(1,0.45791358)%
    \lineheight{1}%
    \setlength\tabcolsep{0pt}%
    \put(0,0){\includegraphics[width=\unitlength,page=1]{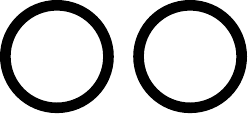}}%
  \end{picture}%
\endgroup%
\ }) = \Delta(  \scalebox{0.15}{%% Creator: Inkscape 1.4.2 (ebf0e940, 2025-05-08), www.inkscape.org
%% Accompanies image file 'unknot+.pdf' (pdf, eps, ps)
\begingroup%
  \makeatletter%
  \providecommand\color[2][]{%
    \errmessage{(Inkscape) Color is used for the text in Inkscape, but the package 'color.sty' is not loaded}%
    \renewcommand\color[2][]{}%
  }%
  \providecommand\transparent[1]{%
    \errmessage{(Inkscape) Transparency is used (non-zero) for the text in Inkscape, but the package 'transparent.sty' is not loaded}%
    \renewcommand\transparent[1]{}%
  }%
  \providecommand\rotatebox[2]{#2}%
  \newcommand*\fsize{\dimexpr\f@size pt\relax}%
  \newcommand*\lineheight[1]{\fontsize{\fsize}{#1\fsize}\selectfont}%
  \ifx\svgwidth\undefined%
    \setlength{\unitlength}{118.31346527bp}%
    \ifx\svgscale\undefined%
      \relax%
    \else%
      \setlength{\unitlength}{\unitlength * \real{\svgscale}}%
    \fi%
  \else%
    \setlength{\unitlength}{\svgwidth}%
  \fi%
  \global\let\svgwidth\undefined%
  \global\let\svgscale\undefined%
  \makeatother%
  \begin{picture}(1,0.46032257)%
    \lineheight{1}%
    \setlength\tabcolsep{0pt}%
    \put(0,0){\includegraphics[width=\unitlength,page=1]{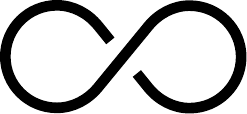}}%
  \end{picture}%
\endgroup%
\ }) - \Delta(\scalebox{0.15}{%% Creator: Inkscape 1.4.2 (ebf0e940, 2025-05-08), www.inkscape.org
%% Accompanies image file 'unknot-.pdf' (pdf, eps, ps)
\begingroup%
  \makeatletter%
  \providecommand\color[2][]{%
    \errmessage{(Inkscape) Color is used for the text in Inkscape, but the package 'color.sty' is not loaded}%
    \renewcommand\color[2][]{}%
  }%
  \providecommand\transparent[1]{%
    \errmessage{(Inkscape) Transparency is used (non-zero) for the text in Inkscape, but the package 'transparent.sty' is not loaded}%
    \renewcommand\transparent[1]{}%
  }%
  \providecommand\rotatebox[2]{#2}%
  \newcommand*\fsize{\dimexpr\f@size pt\relax}%
  \newcommand*\lineheight[1]{\fontsize{\fsize}{#1\fsize}\selectfont}%
  \ifx\svgwidth\undefined%
    \setlength{\unitlength}{116.89614243bp}%
    \ifx\svgscale\undefined%
      \relax%
    \else%
      \setlength{\unitlength}{\unitlength * \real{\svgscale}}%
    \fi%
  \else%
    \setlength{\unitlength}{\svgwidth}%
  \fi%
  \global\let\svgwidth\undefined%
  \global\let\svgscale\undefined%
  \makeatother%
  \begin{picture}(1,0.45377918)%
    \lineheight{1}%
    \setlength\tabcolsep{0pt}%
    \put(0,0){\includegraphics[width=\unitlength,page=1]{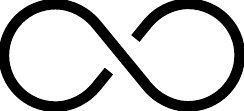}}%
  \end{picture}%
\endgroup%
\ }) = 1-1=0 \ \ \Rightarrow \  \ 
\Delta( \scalebox{0.15}{\ }) =0.$$
For the trefoil $T$, a slightly longer calculation shows that
$$ \Delta(T) = q^{-1} - 1 + q.$$
This implies that $T$ is not isotopic to the unknot.
\begin{figure} \centering
{
\fontsize{9pt}{11pt}\selectfont
   \def\svgwidth{2.5in} 
   %% Creator: Inkscape 1.4.2 (ebf0e940, 2025-05-08), www.inkscape.org
%% Accompanies image file 'orientedskein.pdf' (pdf, eps, ps)
\begingroup%
  \makeatletter%
  \providecommand\color[2][]{%
    \errmessage{(Inkscape) Color is used for the text in Inkscape, but the package 'color.sty' is not loaded}%
    \renewcommand\color[2][]{}%
  }%
  \providecommand\transparent[1]{%
    \errmessage{(Inkscape) Transparency is used (non-zero) for the text in Inkscape, but the package 'transparent.sty' is not loaded}%
    \renewcommand\transparent[1]{}%
  }%
  \providecommand\rotatebox[2]{#2}%
  \newcommand*\fsize{\dimexpr\f@size pt\relax}%
  \newcommand*\lineheight[1]{\fontsize{\fsize}{#1\fsize}\selectfont}%
  \ifx\svgwidth\undefined%
    \setlength{\unitlength}{310.61982883bp}%
    \ifx\svgscale\undefined%
      \relax%
    \else%
      \setlength{\unitlength}{\unitlength * \real{\svgscale}}%
    \fi%
  \else%
    \setlength{\unitlength}{\svgwidth}%
  \fi%
  \global\let\svgwidth\undefined%
  \global\let\svgscale\undefined%
  \makeatother%
  \begin{picture}(1,0.29616998)%
    \lineheight{1}%
    \setlength\tabcolsep{0pt}%
    \put(0,0){\includegraphics[width=\unitlength,page=1]{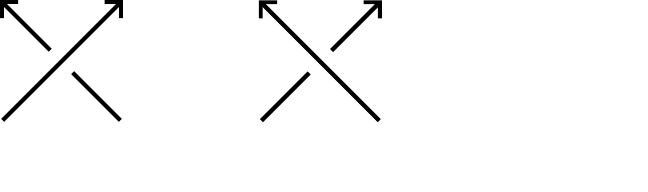}}%
    \put(0.05050422,0.00505656){\color[rgb]{0,0,0}\makebox(0,0)[lt]{\lineheight{1.25}\smash{\begin{tabular}[t]{l}$L_+$\end{tabular}}}}%
    \put(0.4666739,0.00505656){\color[rgb]{0,0,0}\makebox(0,0)[lt]{\lineheight{1.25}\smash{\begin{tabular}[t]{l}$L_-$\end{tabular}}}}%
    \put(0,0){\includegraphics[width=\unitlength,page=2]{orientedskein.pdf}}%
    \put(0.87271921,0.00482905){\color[rgb]{0,0,0}\makebox(0,0)[lt]{\lineheight{1.25}\smash{\begin{tabular}[t]{l}$L_0$\end{tabular}}}}%
  \end{picture}%
\endgroup%

}
\caption{Three links that differ near a crossing}
\label{fig:orientedskein}
\end{figure}

\begin{figure} \centering
{
\fontsize{9pt}{11pt}\selectfont
   \def\svgwidth{1in} 
   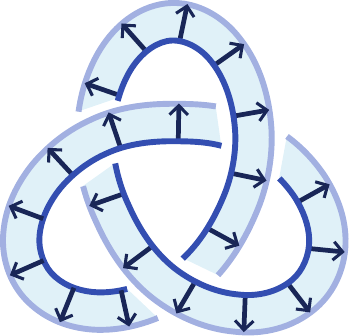
}
\caption{A diagram of the (right-handed) trefoil with its blackboard framing. This corresponds to $\lambda=3$.}
\label{fig:framing}
\end{figure}

Another link invariant is the {\em Jones polynomial} \cite{Jones}, denoted $V(L)$. This was discovered much later than the Alexander polynomial, although it has a very similar characterization. It is normalized by $V(U)=1$ and satisfies the skein relation
$$ q^{-1} V(L_+) - q V(L_-) = (q^{1/2} - q^{-1/2}) V(L_o).$$
The Jones polynomial does not have a simple topological definition, but it admits an alternate skein-theoretic characterization, due to Kauffman \cite{Kauffman}. To describe it, let us first introduce {\em framed links}. A framed link is a link $L$ together with a normal vector field; alternatively, we can push off the link along that field and obtain a parallel copy $L'$ of $L$, and view the framed link as the pair $(L, L')$.  In particular, given a planar diagram $D$ of $L$, we can push off the link along a normal vector inside the plane; this gives the {\em blackboard framing}. See Figure~\ref{fig:framing}.

Up to homotopy, a framing can be specified by a single integer, the linking number $\lambda$ between $L$ and $L'$. We sometimes write the framed link as $\Link = (L, \lambda)$. The framing corresponding to $\lambda =0$ is called the {\em Seifert framing}, and can be different from the blackboard framing. Indeed, the value of $\lambda$ for the blackboard framing is given by the writhe
$$ w=\frac{1}{2} \left(\# \text{positive crossings} - \# \text{negative crossings} \right),$$
where a crossing is called positive or negative according to whether it looks like $L_+$ or $L_-$ in Figure~\ref{fig:orientedskein}.

\begin{figure} \centering
{
\fontsize{9pt}{11pt}\selectfont
   \def\svgwidth{2.2in} 
   %% Creator: Inkscape 1.4.2 (ebf0e940, 2025-05-08), www.inkscape.org
%% Accompanies image file 'unorientedskein.pdf' (pdf, eps, ps)
\begingroup%
  \makeatletter%
  \providecommand\color[2][]{%
    \errmessage{(Inkscape) Color is used for the text in Inkscape, but the package 'color.sty' is not loaded}%
    \renewcommand\color[2][]{}%
  }%
  \providecommand\transparent[1]{%
    \errmessage{(Inkscape) Transparency is used (non-zero) for the text in Inkscape, but the package 'transparent.sty' is not loaded}%
    \renewcommand\transparent[1]{}%
  }%
  \providecommand\rotatebox[2]{#2}%
  \newcommand*\fsize{\dimexpr\f@size pt\relax}%
  \newcommand*\lineheight[1]{\fontsize{\fsize}{#1\fsize}\selectfont}%
  \ifx\svgwidth\undefined%
    \setlength{\unitlength}{304.46337602bp}%
    \ifx\svgscale\undefined%
      \relax%
    \else%
      \setlength{\unitlength}{\unitlength * \real{\svgscale}}%
    \fi%
  \else%
    \setlength{\unitlength}{\svgwidth}%
  \fi%
  \global\let\svgwidth\undefined%
  \global\let\svgscale\undefined%
  \makeatother%
  \begin{picture}(1,0.28243241)%
    \lineheight{1}%
    \setlength\tabcolsep{0pt}%
    \put(0,0){\includegraphics[width=\unitlength,page=1]{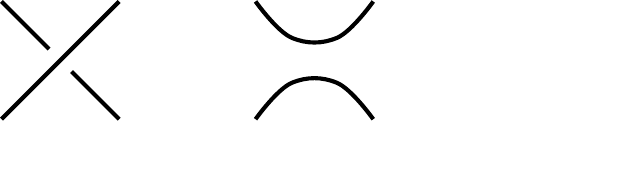}}%
    \put(0.06299362,0.00487595){\color[rgb]{0,0,0}\makebox(0,0)[lt]{\lineheight{1.25}\smash{\begin{tabular}[t]{l}$D$\end{tabular}}}}%
    \put(0.46020586,0.00568394){\color[rgb]{0,0,0}\makebox(0,0)[lt]{\lineheight{1.25}\smash{\begin{tabular}[t]{l}$D_0$\end{tabular}}}}%
    \put(0,0){\includegraphics[width=\unitlength,page=2]{unorientedskein.pdf}}%
    \put(0.86270868,0.00492669){\color[rgb]{0,0,0}\makebox(0,0)[lt]{\lineheight{1.25}\smash{\begin{tabular}[t]{l}$D_1$\end{tabular}}}}%
  \end{picture}%
\endgroup%

}
\caption{A link diagram and its two resolutions at a crossing}
\label{fig:unorientedskein}
\end{figure}

The {\em Kauffman bracket} $\langle D \rangle \in \Z[A, A^{-1}]$ of a planar link diagram $D$ is determined by the properties
\begin{align}
\label{eq:KB0}
 \langle \bigcirc \rangle &= 1,  \\
\label{eq:KB1}
  \langle \bigcirc \sqcup D \rangle &= (-A^2 - A^{-2})\langle D \rangle,
 \end{align}
and a skein relation that, instead of the links in Figure~\ref{fig:orientedskein}, uses those from Figure~\ref{fig:unorientedskein}:
\begin{equation}
\label{eq:KB}
 \langle D \rangle = A\langle D_0 \rangle  + A^{-1} \langle D_1 \rangle.
 \end{equation}

It turns out that the resulting Laurent polynomial $\langle D \rangle$ is not an invariant of the link, but one of the framed link $\Link=(L, w)$, which is represented by the diagram together with its blackboard framing. To get a true link invariant (in fact, the Jones polynomial), we renormalize as follows:
\begin{equation}
\label{eq:JonesKB}
 V(L) = (-A)^{-3w} \langle D \rangle|_{q^{1/2} = A^{-2}}.
  \end{equation}

Both the Alexander and the Jones polynomials can be upgraded to homology theories: knot Floer homology and Khovanov homology, respectively. We leave their definitions for Sections~\ref{sec:hfk} and \ref{sec:Kh}. For now, we mention a few of their properties.

Knot Floer homology, discovered independently by Ozsv\'ath-Szab\'o \cite{Knots} and Rasmussen \cite{RasmussenThesis}, associates to a knot $K$ a bigraded abelian group
\begin{equation}
\label{eq:HFKhat}
\HFKhat(K) = \bigoplus_{i, s \in \Z} \HFKhat_i(K, s).
\end{equation}
The graded Euler characteristic of $\HFKhat$ is the Alexander polynomial:
$$
 \sum_{s, i\in \Z} (-1)^i q^s \cdot \rk_{\Z}  \bigl( \HFKhat_i(K, s) \bigr) = \Delta(K).
$$
 
 Similarly, Khovanov homology \cite{Khovanov} is a bigraded abelian group
\begin{equation}
\label{eq:Kh}
 \Kh(K) = \bigoplus_{i,j} \Kh_{i,j}(K)
 \end{equation}
 whose graded Euler characteristic is the (suitably normalized) Jones polynomial
 $$
 \sum_{i, j\in \Z} (-1)^i q^{j/2} \cdot \rk_{\Z}  \bigl( \Kh_{i, j}(K) \bigr) = (q^{1/2}+q^{-1/2})V(K).
$$

 The knot homologies are more powerful invariants than the corresponding polynomials. For example, a question one can ask about a knot invariant is whether it detects the unknot $U$: if a knot has the same invariant as $U$, is it isotopic to $U$? In the case of the Alexander polynomial, the answer is negative: the Conway knot $C$ shown in Figure~\ref{fig:knots}(e) has $\Delta(C)=1$, just like the unknot. In the case of the Jones polynomial, unknot detection is a famous open problem. However, for the two knot homologies, the question has affirmative answers.
 \begin{theorem}[Ozsv\'ath-Szab\'o \cite{GenusBounds}]
 If a knot $K$ satisfies $\HFKhat(K) \cong \HFKhat(U)$, then $K$ is isotopic to the unknot $U$.
 \end{theorem}
 
 \begin{theorem}[Kronheimer-Mrowka \cite{KMUnknot}]
 If a knot $K$ satisfies $\Kh(K) \cong \Kh(U)$, then $K$ is isotopic to the unknot $U$.
 \end{theorem}

\section{Four-manifolds}
\label{sec:4m}
The classification of manifolds is a major problem in topology. There are two versions of this problem. One can ask about {\em topological manifolds} (spaces that look locally like $\R^n$), in which case the classification should be up to {\em homeomorphism} (bijective continuous map whose inverse is continuous). On the other hand, one can ask about {\em smooth manifolds} (topological manifolds equipped with a smooth structure, i.e., a $C^\infty$ atlas of charts), in which case the classification should be up to {\em diffeomorphism} (bijective $C^\infty$ map whose inverse is $C^\infty$). 

The two classification problems are equivalent for manifolds of dimension up to $3$: any such topological manifold admits a unique smooth structure up to diffeomorphism. In dimensions $4$ and higher, the problems diverge:
\begin{itemize}
\item There exist non-smoothable topological manifolds. The first such example was found by Kervaire in dimension 10 \cite{Kervaire};
\item There exist topological manifolds with several non-diffeomorphic smooth structures. The first example was found by Milnor in dimension 7 \cite{Milnor}.
\end{itemize} 
Although these phenomena exits in higher dimensions, they take a particularly striking form in dimension four. For example:
\begin{itemize}
\item In all dimensions $d \neq 4$, the Euclidean space $\R^d$ admits a unique smooth structure. By contrast, $\R^4$ has uncountably many smooth structures \cite{TaubesPeriodic}; 

\item In all dimensions $d \neq 4$, a compact manifold can admit only finitely many smooth structures \cite{KirbySiebenmann}. By contrast, in dimension $4$ there are compact examples with infinitely many smooth structures \cite{FriedmanMorgan, OkonekV}; 

\item In all dimensions $d \neq 4$, counting the number of smooth structures on the $d$-dimensional sphere $S^d$ has been reduced to a problem in algebraic topology.  By contrast, in dimension $4$ this is a famous open problem, which we now state.
\end{itemize}

\begin{conjecture}[Smooth four-dimensional Poincar\'e Conjecture (SPC4)] Up to diffeomorphism, there exists a unique smooth structure on $S^4$.
\end{conjecture}

The topological four-dimensional Poincar\'e Conjecture (that a four-manifold homotopy equivalent to $S^4$ is actually homeomorphic to $S^4$) was proved by Freedman \cite{Freedman}. In fact, his work produced a complete classification of closed, simply connected topological four-manifolds, up to homeomorphism. (The simply connected restriction is necessary: if we allow the manifolds to have an arbitrary fundamental group $\pi_1$, their classification is impeded by an undecidable problem in group theory, the group isomorphism problem.)

On the smooth side, much progress came from {\em gauge theory}: the study of certain nonlinear, elliptic partial differential equations that originated from physics and exhibit {\em gauge symmetry} (symmetry under the infinite-dimensional group of automorphisms of a bundle). The use of gauge theory in four-dimensional topology dates back to Donaldson \cite{Donaldson}, who looked at self-dual solutions of the Yang-Mills equation. The Seiberg-Witten equations were introduced in 1994 \cite{SeibergWitten1, SeibergWitten2, Witten}, and became an even more effective tool. 

Recall from the introduction that an exotic pair consists of two smooth manifolds that are homeomorphic but not diffeomorphic (in other words, two non-diffeomorphic smooth structures on the same topological manifold). One cannot distinguish the diffeomorphism type by classical topological invariants such as homology or homotopy groups. Rather, in dimension four, the two manifolds in an exotic pair are distinguished roughly as follows. We equip each manifold with a generic Riemannian metric and count the number of solutions to a gauge-theoretic equation. Each solution is counted with a sign, and one shows that (under certain hypotheses) the total count is independent of the metric. It is thus a four-manifold invariant, and if its values are different on the two manifolds, we conclude that the manifolds are not diffeomorphic. 

Of course, to run this strategy one needs a good source of examples for which the invariants are computable. One such source is algebraic geometry.  Projective algebraic surfaces provide many examples of four-manifolds. The simplest ones are the complex projective space $\CP^2$ and its blow-ups. Topologically, these can be expressed as connected sums 
$$ \CP^2 \# n \bCP,$$
where $\bCP$ is $\CP^2$ with the orientation reversed. Another famous example is the $K3$ surface
$$ K3 =\{[z_0: z_1: z_2: z_3] \in \CP^3 \mid z_0^4 + z_1^4+z_2^4 + z_3^4=0\}.$$

Solutions to the Yang-Mills equations on projective algebraic surfaces are related to holomorphic vector bundles. This enables the computation of the Donaldson invariants (signed counts of solutions). Similarly, solutions to the Seiberg-Witten equations on projective algebraic surfaces are related to divisors, and this helps the computation of the corresponding Seiberg-Witten invariants. As sample applications, one can prove the existence of exotic smooth structures on $\CP^2 \# 9 \bCP$ or on $K3$; see \cite{Donaldson2, OkonekV, FriedmanMorgan, FS96}. Seiberg-Witten theory also has deep connections to symplectic geometry \cite{TaubesSymp, TaubesGW}. By combining these methods with cut-and-paste techniques, one can construct exotic smooth structures on $\CP^2 \# k \bCP$ for all $k \geq 2$; see \cite{AP}.

As a rule of thumb, the larger the four-manifold (in terms of the size of its homology), the easier it is to find exotic smooth structures. On relatively simple four-manifolds such as $$S^4, S^2 \times S^2, \CP^2, \CP^2 \# \bCP, S^1 \times S^3, T^4,$$ the existence of exotic smooth structures is still an open problem. 

Another key factor in the construction of exotica is the type of intersection form. On a closed, simply-connected four-manifold $X$, there is a non-degenerate bilinear form on second homology
$$ H_2(X; \Z) \times H_2(X; \Z) \to \Z,$$
given by the intersection product. All known exotic examples of simply-connected four-manifolds have the property that this form is indefinite. When the form is definite, we might as well focus on the positive-definite case (by changing the orientation if necessary). It then follows from Freedman's work that the manifold is homeomorphic to a connected sum $\#^n\CP^2$. The following question remains unsolved.

\begin{question}
\label{q:definite}
Does there exist an exotic smooth structure on $\#^n\CP^2$ for some $n \geq 0$?
\end{question}

Note that SPC4 is a particular case of this problem, corresponding to $n=0$.

\section{Kirby diagrams}
\label{sec:kirby}
The connection between knot theory and four-manifolds runs through three-manifolds, so let us first describe how one can construct three-manifolds using links. Let $\Link = (L, \lambda)$ be a framed link in $S^3$ with $\ell$ components, and let $\nu(L)$ a standard tubular neighborhood of $L$. The result of {\em surgery} on $\Link$  is the three-manifold
$$ S^3(\Link) = (S^3 \setminus \nu(L)) \cup_{\del \nu(L)} \bigsqcup_{i=1}^{\ell} (S^1 \times D^2),$$
obtained from the  complement of $\nu(L)$ by gluing back $\ell$ solid tori $S^1 \times D^2$ along their boundaries. The gluing is done so that the meridians $* \times \del D^2$ are attached to the components of the parallel copy of $L$ specified by the framing.

\begin{theorem}[Lickorish-Wallace \cite{lickorish, Wallace}]
\label{thm:LW}
Any closed, oriented three-manifold is the result of surgery on a framed link in $S^3$.  
\end{theorem}

\begin{figure} \centering
{
\fontsize{9pt}{11pt}\selectfont
   \def\svgwidth{4in} 
   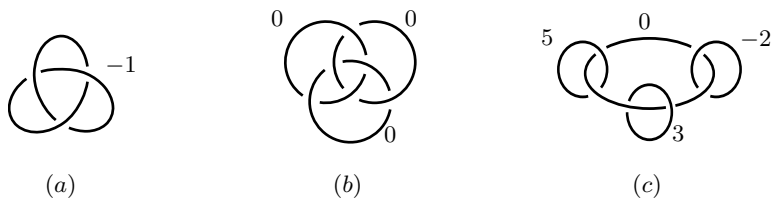
}
\caption{Surgery diagrams for: (a) the Poincar\'e homology sphere; (b) the $3$-torus $T^3=S^1 \times S^1 \times S^1$; (c) the Poincar\'e sphere again}
\label{fig:surgery}
\end{figure}

Pictorially, we represent the three-manifold by a link with integers attached to each component; these determine the framing as noted in Section~\ref{sec:knots}. The result is a {\em surgery diagram} for the manifold. Some examples are shown in Figure~\ref{fig:surgery}. The surgery diagram is not unique; many framed links represent the same three-manifold. 

\begin{figure} \centering
{
\fontsize{9pt}{11pt}\selectfont
   \def\svgwidth{2.1in} 
   %% Creator: Inkscape 1.3.2 (091e20e, 2023-11-25), www.inkscape.org
%% Accompanies image file 'trace.pdf' (pdf, eps, ps)
\begingroup%
  \makeatletter%
  \providecommand\color[2][]{%
    \errmessage{(Inkscape) Color is used for the text in Inkscape, but the package 'color.sty' is not loaded}%
    \renewcommand\color[2][]{}%
  }%
  \providecommand\transparent[1]{%
    \errmessage{(Inkscape) Transparency is used (non-zero) for the text in Inkscape, but the package 'transparent.sty' is not loaded}%
    \renewcommand\transparent[1]{}%
  }%
  \providecommand\rotatebox[2]{#2}%
  \newcommand*\fsize{\dimexpr\f@size pt\relax}%
  \newcommand*\lineheight[1]{\fontsize{\fsize}{#1\fsize}\selectfont}%
  \ifx\svgwidth\undefined%
    \setlength{\unitlength}{324.21843762bp}%
    \ifx\svgscale\undefined%
      \relax%
    \else%
      \setlength{\unitlength}{\unitlength * \real{\svgscale}}%
    \fi%
  \else%
    \setlength{\unitlength}{\svgwidth}%
  \fi%
  \global\let\svgwidth\undefined%
  \global\let\svgscale\undefined%
  \makeatother%
  \begin{picture}(1,0.76735605)%
    \lineheight{1}%
    \setlength\tabcolsep{0pt}%
    \put(0,0){\includegraphics[width=\unitlength,page=1]{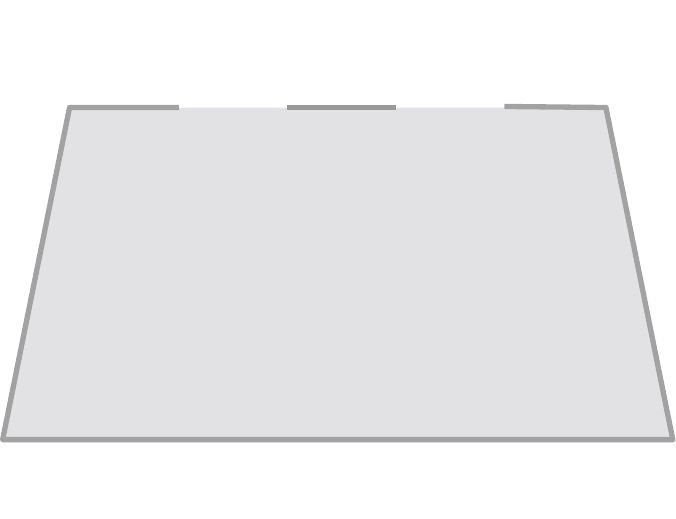}}%
    \put(0.7265132,0.21461178){\color[rgb]{0.20392157,0.70588235,0.65490196}\makebox(0,0)[lt]{\lineheight{1.25}\smash{\begin{tabular}[t]{l}$\nu(L)$\end{tabular}}}}%
    \put(0,0){\includegraphics[width=\unitlength,page=2]{trace.pdf}}%
    \put(0.08945693,0.003109){\color[rgb]{0,0,0}\makebox(0,0)[lt]{\lineheight{1.25}\smash{\begin{tabular}[t]{l}$D^4$\end{tabular}}}}%
    \put(0.0820094,0.17427468){\color[rgb]{0,0,0}\makebox(0,0)[lt]{\lineheight{1.25}\smash{\begin{tabular}[t]{l}$S^3$\end{tabular}}}}%
    \put(0,0){\includegraphics[width=\unitlength,page=3]{trace.pdf}}%
    \put(0.68999609,0.73804249){\color[rgb]{0,0,0}\makebox(0,0)[lt]{\lineheight{1.25}\smash{\begin{tabular}[t]{l}$2$-handle\end{tabular}}}}%
    \put(0,0){\includegraphics[width=\unitlength,page=4]{trace.pdf}}%
  \end{picture}%
\endgroup%

}
\caption{The trace of a surgery. The disk $D^4$ is depicted as the half-space below the indicated plane.}
\label{fig:trace}
\end{figure}

There is also a four-manifold (with boundary) associated to surgery on a framed link $\Link$. This is called the {\em trace of the surgery}, denoted $X(\Link)$, and is constructed as follows:
$$ X(\Link) = D^4 \cup_{\del \nu(L)} \bigsqcup_{i=1}^{\ell} (D^2 \times D^2),$$
Here, we view $\nu(L)$ as part of the boundary $\del D^4 = S^3$, and each copy of $D^2 \times D^2$ (called a {\em $2$-handle}) is attached along $\del D^2 \times D^2=S^1 \times D^2$ to a component of $\nu(L)$. (The identification of such a component with $S^1 \times D^2$ is determined by the framing.) See Figure~\ref{fig:trace}.

The result of attaching the $2$-handle to $D^4$ is {\em a priori} a four-manifold with corners. In the definition of $X(\Link)$ we implicitly smooth the corners. The result is a compact four-manifold with boundary, and one can check that $\del X(\Link) = S^3(\Link)$. In light of this, Theorem~\ref{thm:LW} has the following consequence.

\begin{corollary}
\label{cor:3bdry}
Every closed, oriented three-manifold is the boundary of a compact, oriented four-manifold.
\end{corollary}

\begin{figure} \centering
{
\fontsize{9pt}{11pt}\selectfont
   \def\svgwidth{5in} 
  %% Creator: Inkscape 1.3.2 (091e20e, 2023-11-25), www.inkscape.org
%% Accompanies image file 'morse.pdf' (pdf, eps, ps)
\begingroup%
  \makeatletter%
  \providecommand\color[2][]{%
    \errmessage{(Inkscape) Color is used for the text in Inkscape, but the package 'color.sty' is not loaded}%
    \renewcommand\color[2][]{}%
  }%
  \providecommand\transparent[1]{%
    \errmessage{(Inkscape) Transparency is used (non-zero) for the text in Inkscape, but the package 'transparent.sty' is not loaded}%
    \renewcommand\transparent[1]{}%
  }%
  \providecommand\rotatebox[2]{#2}%
  \newcommand*\fsize{\dimexpr\f@size pt\relax}%
  \newcommand*\lineheight[1]{\fontsize{\fsize}{#1\fsize}\selectfont}%
  \ifx\svgwidth\undefined%
    \setlength{\unitlength}{477.75708777bp}%
    \ifx\svgscale\undefined%
      \relax%
    \else%
      \setlength{\unitlength}{\unitlength * \real{\svgscale}}%
    \fi%
  \else%
    \setlength{\unitlength}{\svgwidth}%
  \fi%
  \global\let\svgwidth\undefined%
  \global\let\svgscale\undefined%
  \makeatother%
  \begin{picture}(1,0.37893993)%
    \lineheight{1}%
    \setlength\tabcolsep{0pt}%
    \put(0,0){\includegraphics[width=\unitlength,page=1]{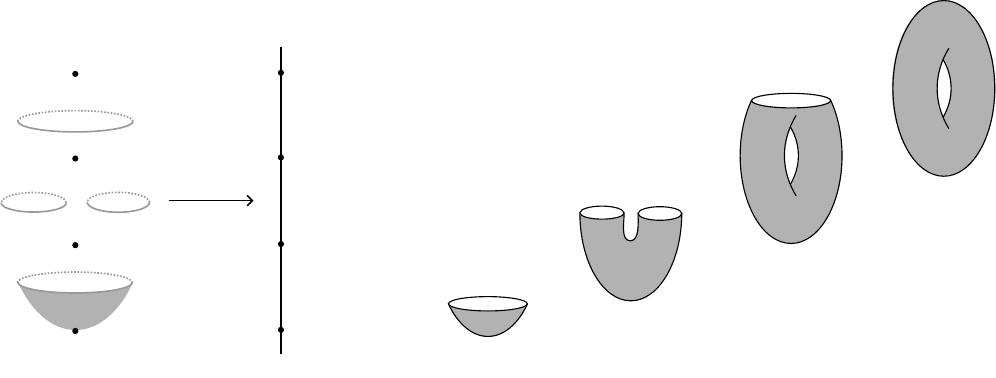}}%
    \put(0.33565134,0.00378968){\color[rgb]{0,0,0}\makebox(0,0)[lt]{\lineheight{1.25}\smash{\begin{tabular}[t]{l}$\emptyset$\end{tabular}}}}%
    \put(0,0){\includegraphics[width=\unitlength,page=2]{morse.pdf}}%
    \put(0.66605522,0.07242848){\color[rgb]{0,0,0}\makebox(0,0)[lt]{\lineheight{1.25}\smash{\begin{tabular}[t]{l}$=$\end{tabular}}}}%
    \put(0,0){\includegraphics[width=\unitlength,page=3]{morse.pdf}}%
  \end{picture}%
\endgroup%

}
\caption{A handle decomposition of the $2$-torus from a Morse function. As we pass critical points from bottom to top, the sublevel set changes by attaching first a $0$-handle, then a $1$-handle, then another $1$-handle, and finally a $2$-handle.}
\label{fig:morse}
\end{figure}

Surgery traces form a key part of {\em Kirby diagrams}, the standard pictorial representations of four-manifolds. To introduce Kirby diagrams, let us first review a few notions of Morse theory. On a closed, connected manifold $X$ of some dimension $n \geq 0$, a generic smooth function $f: X \to \R$ is {\em Morse}, that is, all of its critical points are non-degenerate. As we follow the sublevel sets $\{f(x) \leq C$ as $C$ goes from $-\infty$ to $\infty$, their diffeomorphism type changes only as we pass critical values, and in that case it changes by attaching handles. Specifically, when we encounter a critical point of index $k$, we attach a $k$-handle $D^k \times D^{n-k}$ along $\del D^k \times D^{n-k}$. One can re-arrange the handles so that they appear monotonically ($k$-handles before $l$-handles for $k < l$), and we have a single $0$-handle and a single $n$-handle. We refer to \cite{Morse} and \cite{MilnorHcob} for more details. See Figure~\ref{fig:morse} for the example of the $2$-torus.

In the case of a closed four-manifold $X$, we get a decomposition into handles as follows:
\begin{itemize}
\item one $0$-handle, which is a copy of $D^4$;
\item several $1$-handles, which produce a boundary connected sum $\natural^k(S^1 \times D^3)$, whose boundary is $\#^k(S^1 \times S^2)$;
\item several $2$-handles, which are attached along the neighborhood of a framed link $\Link \subset \#^k(S^1 \times S^2)$. When $k=0$, this corresponds to the surgery trace discussed earlier;
\item some $3$-handles, which are attached along neighborhoods of $2$-spheres in the boundary of the manifold obtained at the previous step. We assume that, after attaching the $3$-handles, the new boundary is  $S^3$;
\item a $4$-handle, which is a copy of $D^4$ attached along its whole boundary.
\end{itemize}

\begin{figure} \centering
{
\fontsize{9pt}{11pt}\selectfont
   \def\svgwidth{2.5in} 
   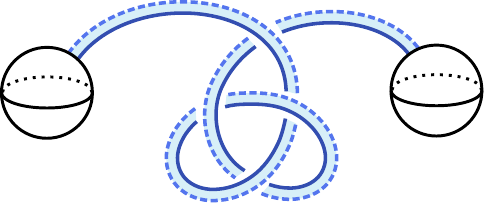
}
\caption{A Kirby diagram (representing $S^4$). The two spheres are identified by reflection about the vertical plane.}
\label{fig:kirby1}
\end{figure}

The Kirby diagram associated to this decomposition is simply a picture of the framed link $\Link \subset \#^k(S^1 \times S^2)$. To draw $\#^k(S^1 \times S^2)$, we represent $S^3$ as $\R^3 \cup \{ \infty\}$, then remove $k$ pairs of balls from $\R^3$ and identify the boundaries of the balls in each pair. The framing is indicated by a dashed parallel copy of the link, or (when $k=0$) by marking an integer next to each link component. An example of a Kirby diagram is shown in Figure~\ref{fig:kirby1}.

Interestingly, under the assumption that $3$-handles can be attached so that the new boundary is $S^3$, the way we attach these handles is essentially unique, and does not require extra data in the diagram. Indeed, the boundary after we attach the $2$-handles must be some connected sum $\#^l(S^1 \times S^2)$, and the uniqueness boils down to a theorem of Laudenbach and Po{\'e}naru \cite{LaudenbachPoenaru}, which says that every diffeomorphism of $\#^l(S^1 \times S^2)$ extends to the handlebody $\natural^l(S^1 \times D^3)$.

Of course, the presentation of a four-manifold by a Kirby diagram is not unique. It can be shown that any two such presentations are related by a sequence of certain moves \cite{KirbyCalc}. This is the subject of {\em Kirby calculus}; see \cite{GompfStipsicz} for an introduction.

\begin{figure} \centering
{
\fontsize{9pt}{11pt}\selectfont
  \def\svgwidth{5in} 
   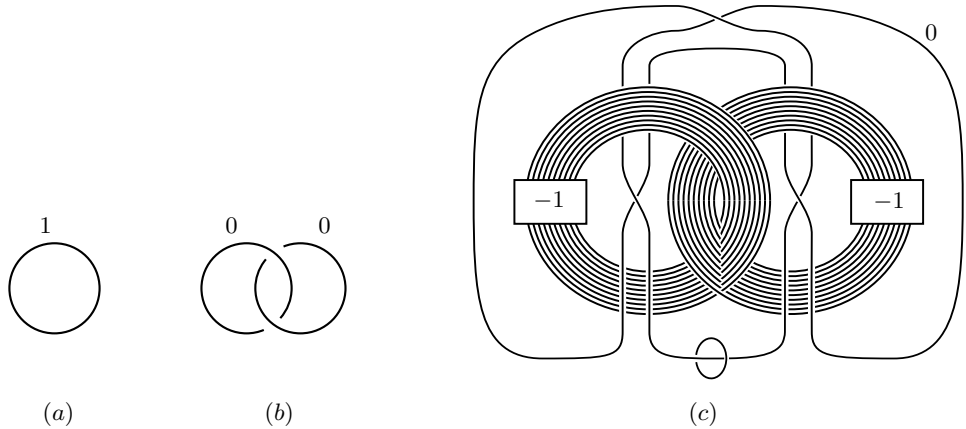

}
\caption{Kirby diagrams for: (a) The complex projective space $\CP^2$; (b) $S^2 \times S^2$; (c) The $K3$ surface. The third picture is based on Figure 2 in \cite{MMP}; all the circles there have framing $-2$, except for the trefoil with framing $0$. The boxes marked $-1$ indicate a negative full twist of all the strands in the box. }
\label{fig:kirby2}
\end{figure}

In this paper we will mostly deal with the case when the handle decomposition has no $1$-handles, so the link $\Link$ is in $S^3$. Such a manifold is called {\em geometrically simply connected}, and it is an open problem whether every closed, simply connected four-manifold has this property. In practice, many simply connected manifolds have it. Figure~\ref{fig:kirby2} shows Kirby diagrams (without $1$- or $3$-handles) for $\CP^2$, $S^2 \times S^2$, and the $K3$ surface. All of these manifolds are obtained from a surgery trace $X(\Link)$ with boundary $S^3(\Link)=S^3$ by attaching a copy of $D^4$.

There is also a variant of Kirby diagrams for compact, connected four-manifolds $X$ with non-empty boundary. A handle decomposition for such a manifold can be taken to have a single $0$-handle, some $1$-, $2$-, and $3$-handles, and no $4$-handles. We again draw it as a framed link $\Link \subset \#^k(S^1 \times S^2)$. In this case we no longer have the assumption on the boundary after attaching $3$-handles; this boundary can be any three-manifold $Y=\del X$. We can no longer apply the Laudenbach-Po{\'e}naru theorem so the way the $3$-handles are attached has to be specified in the text. 

In particular, if we draw a framed link $\Link$ in $S^3$, this can be viewed as a Kirby diagram for either:
\begin{enumerate}
\item[(a)] a closed four-manifold, in case $S^3(\Link) = \#^l(S^1 \times S^2)$ for some $l$; or
\item[(b)] a compact four-manifold with boundary $S^3(\Link)$, in general.
\end{enumerate}

\section{Heegaard Floer theory}
\label{sec:HF}

\subsection{Floer homologies.}
\label{sec:3mi}
Recall from Section~\ref{sec:4m} that counting solutions to the Yang-Mills or Seiberg-Witten equations produces the Donaldson and Seiberg-Witten invariants, respectively, which can detect exotic smooth structures.     To compute these invariants, in some cases one can use connections to geometry, but a more widely applicable method is through cut-and-paste techniques. The idea is to break the four-manifold into simpler pieces, compute invariants of these pieces, and investigate how these behave under gluing.

The pieces are four-manifolds with boundary, and their invariants will not be numerical (as for closed manifolds), but rather elements in a group associated to the boundary, called {\em Floer homology}. This group is an invariant of three-manifolds constructed roughly as follows. Given a three-manifold $Y$, we look at  solutions to the relevant equations (Yang-Mills or Seiberg-Witten) on the cylinder $\R \times Y$ that are invariant under translation by $\R$. We define a chain complex whose generators are these solutions, and whose differential is of the form 
$$ \del x = \sum_y n(x, y) y,$$
with $n(x, y) \in \Z$ being the count of (not necessarily translation invariant) solutions on $\R \times Y$ that limit to $x$ and $y$ as we go to $-\infty$ and $+\infty$ in the $\R$ direction. The homology of the complex is called Floer homology. In the Yang-Mills context, this was introduced in \cite{FloerInstanton} and in the Seiberg-Witten context in \cite{MarcolliWang, Spectrum, FroyshovSW, KMBook}. 

Let us denote the Floer homology theory of $Y$ (in some setting) by $F(Y)$. Then, the invariant of a four-manifold $W$ with boundary $Y$ is 
$$F(W)=[\sum_x m(x) x] \in F(Y),$$
where $m(x)$ is the count of solutions on the manifold with cylindrical end $X \cup ([0, \infty) \times Y)$, converging to $x$ as we go to the infinite end.  We can alternatively think of $F(W)$ as a homomorphism 
$$ F(\emptyset) = \Z \to F(Y)$$
sending $1$ to the count above. More generally, if we have a cobordism $W$ between three-manifolds $Y$ and $Y'$ (that is, a four-manifold $W$ with $\del W = (-Y) \sqcup Y'$), we get a map
$$ F(W): F(Y) \to F(Y').$$
This assignment is functorial: composition of cobordisms results in composition of the maps. We obtain a structure  called a {\em topological quantum field theory} (TQFT). 

If a closed four-manifold $W$ is decomposed into pieces as
$$ W = W_1 \cup_{Y_1} W_2 \cup_{Y_2} \cdots \cup_{Y_{n-1}} W_n,$$
its invariant is given by
\begin{equation}
\label{eq:FW}
 F(W) = F(W_n) \circ \cdots \circ F(W_2) \circ F(W_1).
 \end{equation}
Thus, it suffices to understand the invariants of the pieces.

\subsection{Heegaard Floer homology}
The first step in making use of Floer-theoretic cut-and-paste techniques is to understand the Floer homology of three-manifolds. This is non-trivial, as it still involves solving PDEs on $\R \times Y$. It is useful to make an additional cut in the three-manifold: split it along a surface $\Sigma$ into two handlebodies $U_0$ and $U_1$. This is called a {\em Heegaard splitting}. By stretching the metric along  $\Sigma$, the original equations on $\R \times Y$ become maps
\begin{equation}
\label{eq:u}
 u: \R \times [0, 1] \to M(\Sigma),
 \end{equation}
where $M(\Sigma)$ is the moduli space of dimensionally-reduced solutions on $\Sigma$; that is, solutions to the PDEs on $\R^2 \times \Sigma$ that are invariant under translation in both $\R$ directions. Furthermore, it turns out that $M(\Sigma)$ has a natural symplectic structure, and the maps in \eqref{eq:u} need to be pseudo-holomorphic: They satisfy a nonlinear analogue of the Cauchy-Riemann equation, with respect to an almost complex structure compatible with the symplectic form. 

This heuristic shows that the gauge-theoretic Floer homology of $Y$ should be isomorphic to a symplectic invariant called {\em Lagrangian Floer homology}. This proposal was made in the Yang-Mills setting in \cite{AtiyahFloer}, and came to be known as the Atiyah-Floer conjecture. It was even more fruitful in the Seiberg-Witten setting, where Ozsv\'ath and Szab\'o developed {\em Heegaard Floer theory} \cite{HolDisk, HolDiskTwo}.

The dimensionally-reduced Seiberg-Witten equations on a surface $\Sigma$ are called the {\em vortex equations}, and their moduli spaces are symmetric products of the surface \cite{JaffeTaubes, Bradlow}. Inspired by this, Ozsv\'ath and Szab\'o defined the Heegaard Floer homology of a three-manifold $Y$ with a Heegaard splitting $U_0 \cup_{\Sigma} U_1$ by counting pseudo-holomorphic strips as in \eqref{eq:u}, with 
$$M(\Sigma) =\Sym^g(\Sigma) = (\Sigma \times \cdots \times \Sigma)/S_g.$$
Here, $g$ is the genus of $\Sigma$; on the right hand side we take the product of $g$ copies of $\Sigma$ and divide it by the action of the symmetric group $S_g$. The handlebody $U_0$ can be specified by a collection of $g$ simple closed curves $\alpha_1, \dots, \alpha_g$ on $\Sigma$, such that attaching disks to $\Sigma$ along these curves and then filling in with a three-ball $D^3$ results in $U_0$. Similarly, the handlebody $U_1$ can be specified by a collection of curves $\beta_1, \dots, \beta_g$. From here we construct Lagrangian submanifolds
$$ \Ta = \alpha_1 \times \cdots \times \alpha_g, \ \ \Tb = \beta_1 \times \cdots \times \beta_g \subset \Sym^g(\Sigma).$$

The generators of the Heegaard Floer complex are intersection points $\x \in \Ta \cap \Tb$, and the differential counts pseudo-holomorphic strips in $M(\Sigma)$ with boundaries on $\Ta$ and $\Tb$. The homology of the complex is Heegaard Floer homology, $\HF(Y)$. (In \cite{HolDisk} there are actually several flavors of this construction, denoted $\HFhat$, $\HFp$, $\HFm$, $\HFinfty$, which differ in how they keep track of a basepoint $z \in \Sigma$. We will not get into the details here, and write any of these versions as $\HF$.)

Compared to Seiberg-Witten theory, Heegaard Floer homology has the advantage that at least the generators of the chain complex are easy to understand topologically. They correspond to $g$-tuples of points on $\Sigma$:
$$\x = \{x_1, \dots, x_g\}$$
with $x_i \in \alpha_i \cap \beta_{\sigma(g)}$ for some permutation $\sigma$ of $\{1, \dots, g\}$. Computing the differential is still challenging, because it involves solving the nonlinear Cauchy-Riemann equations.

In \cite{HolDiskFour}, Ozsv\'ath and Szab\'o constructed invariants of four-manifolds (with or without boundary) that fit into the TQFT framework sketched in Section~\ref{sec:3mi}. They are defined by starting from a handle decomposition of the four-manifold, as in Section~\ref{sec:kirby}, associating maps to each handle, and composing these maps as in \eqref{eq:FW}. This offers a first indication that the perspective of viewing four-manifolds in terms of Kirby diagrams may be useful for computing invariants.

We mention that Heegaard Floer homology was later proved to be isomorphic to Seiberg-Witten (monopole) Floer homology; see \cite{KLT, CGH}. Similarly, the Ozsv\'ath-Szab\'o invariants of four-manifolds are expected to be equal to the Seiberg-Witten invariants. This remains a conjecture, but for most purposes the Ozsv\'ath-Szab\'o invariants are a perfectly good substitute for Seiberg-Witten: They have been proved to have many of the same properties, and calculations for specific classes of four-manifolds have yielded the same results. See for example \cite{OSSymp, JabukaMark}.

\subsection{Knot Floer homology.}
\label{sec:hfk}
In Section~\ref{sec:knots} we mentioned an invariant of knots called {\em knot Floer homology}, denoted $\HFKhat$. Its original construction in \cite{Knots} and \cite{RasmussenThesis}, is a variant of that of Heegaard Floer homology. We start with a Heegaard splitting for $S^3$, with $\alpha$ and $\beta$ curves as before, and we encode the knot by marking two basepoints on the surface, away from the curves. Specifically, recall that each handlebody is obtained from the surface $\Sigma$ by adding some disks and a three-ball; in each three-ball we draw a segment joining the basepoints, and the knot is the union of these two segments. 

Knots are simpler objects than three-manifolds, so it should not be surprising that knot Floer homology is (in many ways) easier to study than Heegaard Floer homology. In fact, a completely combinatorial description of $\HFKhat$ was given in \cite{MOS}, based on the idea of generalizing the original set-up by allowing more basepoints and curves. This allows for a judicious choice of knot presentations called grid diagrams.

\begin{figure} \centering
{
\fontsize{9pt}{11pt}\selectfont
   \def\svgwidth{3.7in} 
   %% Creator: Inkscape 1.4.2 (ebf0e940, 2025-05-08), www.inkscape.org
%% Accompanies image file 'grid.pdf' (pdf, eps, ps)
\begingroup%
  \makeatletter%
  \providecommand\color[2][]{%
    \errmessage{(Inkscape) Color is used for the text in Inkscape, but the package 'color.sty' is not loaded}%
    \renewcommand\color[2][]{}%
  }%
  \providecommand\transparent[1]{%
    \errmessage{(Inkscape) Transparency is used (non-zero) for the text in Inkscape, but the package 'transparent.sty' is not loaded}%
    \renewcommand\transparent[1]{}%
  }%
  \providecommand\rotatebox[2]{#2}%
  \newcommand*\fsize{\dimexpr\f@size pt\relax}%
  \newcommand*\lineheight[1]{\fontsize{\fsize}{#1\fsize}\selectfont}%
  \ifx\svgwidth\undefined%
    \setlength{\unitlength}{587.09250623bp}%
    \ifx\svgscale\undefined%
      \relax%
    \else%
      \setlength{\unitlength}{\unitlength * \real{\svgscale}}%
    \fi%
  \else%
    \setlength{\unitlength}{\svgwidth}%
  \fi%
  \global\let\svgwidth\undefined%
  \global\let\svgscale\undefined%
  \makeatother%
  \begin{picture}(1,0.43494084)%
    \lineheight{1}%
    \setlength\tabcolsep{0pt}%
    \put(0,0){\includegraphics[width=\unitlength,page=1]{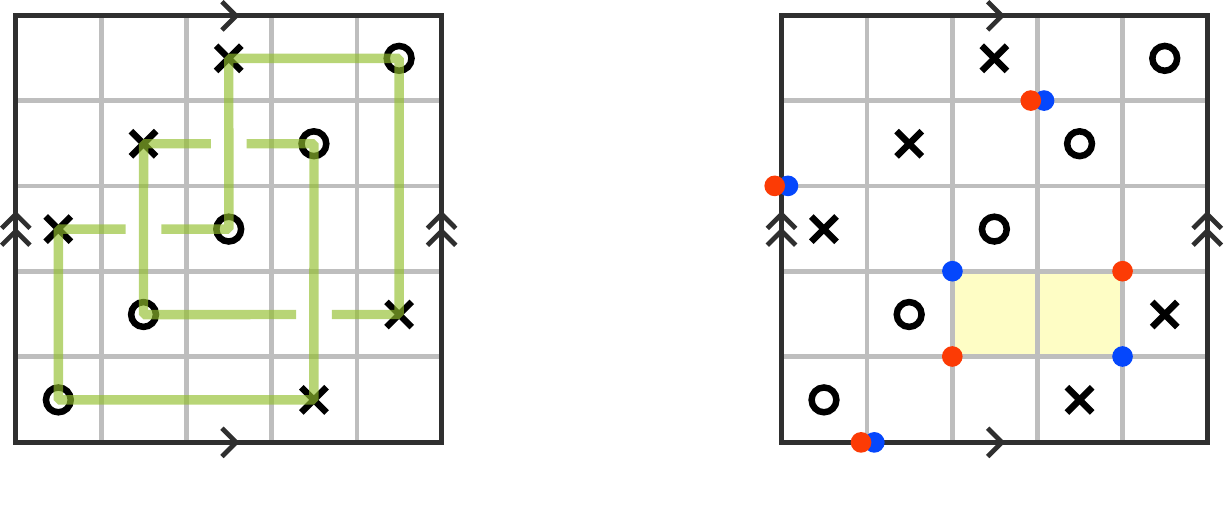}}%
    \put(0.17253527,0.00308391){\color[rgb]{0,0,0}\makebox(0,0)[lt]{\lineheight{1.25}\smash{\begin{tabular}[t]{l}$(a)$\end{tabular}}}}%
    \put(0.80473269,0.00308391){\color[rgb]{0,0,0}\makebox(0,0)[lt]{\lineheight{1.25}\smash{\begin{tabular}[t]{l}$(b)$\end{tabular}}}}%
  \end{picture}%
\endgroup%

}
\caption{(a) A grid diagram representing the trefoil. (b) An empty rectangle from $\x$ to $\y$ on this grid. The generator $\x$ is the collection of red dots, and the generator $\y$ is the collection of blue dots.}
\label{fig:grid}
\end{figure}

A (toroidal) {\em grid diagram} is an $n \times n$ grid on the surface of a torus (viewed as a square with the opposite sides identified), with the following data:
\begin{itemize}
\item a collection of $n$ parallel horizontal curves $\alpha_1, \dots, \alpha_n$, splitting the torus into $n$ rows;
\item a collection of $n$ parallel horizontal curves $\beta_1, \dots, \beta_n$, splitting the torus into $n$ columns;
\item a collection of $n$ markings $O_1, \dots, O_n$, such that each row and each column contain exactly one of these markings;
\item another collection of $n$ markings $X_1, \dots, X_n$, such that each row and each column contain exactly one of these markings.
\end{itemize}
Let us join the $O$ to the $X$ markings by vertical segments inside each column, and by horizontal segments inside each row. When the segments cross each other, we draw the vertical ones on top. The result is a planar diagram for a link $L$. See Figure~\ref{fig:grid}(a) for an example.

We now form a chain complex freely generated over $\Z$ by $n$-tuples of points
$$\x = \{x_1, \dots, x_n\}$$
with $x_i \in \alpha_i \cap \beta_{\sigma(n)}$ for some permutation $\sigma$. Note that there are exactly $n!$ generators. (These come with two gradings, which eventually produce the bi-grading on $\HFKhat$ from \eqref{eq:HFKhat}; we will not give their definitions here.)

The differential on the complex is
$$\del  \x = \sum_{\y} r(\x, \y) \y,$$
where $r(\x, \y)$ is the count of empty rectangles between $\x$ and $\y$. Such rectangles can exist only when the coordinates of $\x$ and $\y$ differ in exactly two rows (and in exactly two columns). Then, those two rows and two columns split the torus into four rectangles. By an orientation convention, the two where the rectangle is on the left of the segment drawn from $\x$ to $\y$ on each row are considered to go ``from $\x$ to $\y$.'' If such a rectangle does not contain any $O$- or $X$-markings, nor any other coordinates of $\x$ and $\y$, it is counted (with a certain sign) in $r(\x, \y)$. See Figure~\ref{fig:grid}(b) for an example.

The homology of this complex is denoted $\HFKt$, and is not quite an invariant of the link. Nevertheless, one can show that if the link has $\ell$ components, then $\HFKt$ is obtained from $\HFKhat$ by tensoring with $n-\ell$ copies of a free abelian group of rank two.  

There are other versions of this construction, where one counts rectangles that may contain $O$- or $X$-markings. These give rise to different flavors of knot Floer homology (sometimes called {\em grid homology}).  A comprehensive introduction is the textbook \cite{OSSgrids}.

The perspective of grid homology yields the following result.

\begin{theorem}[\cite{MOS}]
\label{thm:HFK}
Knot Floer homology (in all its versions) is algorithmically computable. 
\end{theorem}

Indeed, there exist computer programs that use the grid diagram description to calculate $\HFKhat$ \cite{BaldwinGillam}. The key fact that enabled this description was that isolated holomorphic disks in $\Sym^n(T^2)$ with boundary conditions on $\Ta$ and $\Tb$ are in one-to-one correspondence to empty rectangles. In broad strokes, we can say that working with knots and grids resulted in a setting where gauge-theoretic equations can be solved explicitly: each empty rectangle is a pseudo-holomorphic strip in the moduli space of vortices, and hence corresponds to a solution of the Seiberg-Witten equations on $\R^2 \times T^2$.

\subsection{Surgery formulas.} 
\label{sec:surgeryHF}
We now arrive at the crux of the matter: Once we have a link invariant that we understand reasonably well (knot Floer homology), how can this help us compute a four-manifold invariant (in this case, the Ozsv\'ath-Szab\'o invariant, which is a replacement for the Seiberg-Witten invariant)? The connection goes through three-manifold  invariants (Heegaard Floer homology), using what are called surgery formulas. 

In Section~\ref{sec:kirby} we mentioned that every three-manifold $Y$ can be expressed as surgery on a framed link $ \Link \subset S^3$. In the case where $\Link=(L, \lambda)$ is a knot, the Heegaard Floer homology of $Y=S^3(\Link)$ was related to the knot Floer homology of $L$ by the integral surgery formula proved by Ozsv\'ath and Szab\'o in \cite{IntSurg}. This says that $\HF(S^3(\Link))$ is the homology of a mapping cone
$$ A(L) \to A(\emptyset),$$
where $A$ is a certain version of knot Floer homology, and the framing $\lambda$ is involved in defining the map itself. Moreover, the formula identifies the four-manifold invariant of the trace $X(\Link)$ with the class of a given element of $A(\emptyset)$ in the mapping cone above.

This surgery formula was generalized to links by the author and Ozsv\'ath in \cite{LinKSurg}. In that case, instead of a mapping cone we have a mapping hypercube (an iterated mapping cone) involving the knot Floer homologies of the link and all its sublinks, as well as maps, chain homotopies and higher chain homotopies connecting them. If $\Link' \subset \Link$ is a sublink, then the Heegaard Floer map induced by the cobordism from $S^3(\Link')$ to $S^3(\Link)$ (given by surgery on $\Link - \Link'$) is identified with the inclusion of a certain subcomplex of the hypercube.

Theorem~\ref{thm:HFK} provides a combinatorial description of the groups that sit at the vertices of the mapping hypercube mentioned above. There is more work needed to describe the maps and (higher) homotopies, but this was accomplished in \cite{MOT}, in terms of counting other geometric shapes on the grid diagram. The result was a combinatorial description of the Heegaard Floer homology of $S^3(\Link)$, and of related cobordism maps. Thus, if we have a four-manifold that is made only of two-handles (that is, the trace of a link surgery $X(\Link)$), we obtain an expression for its Ozsv\'ath-Szab\'o invariant. This kind of expression can be generalized to all four-manifolds, using handle decompositions and the composition rule \eqref{eq:FW}. We have the following consequence. 

\begin{theorem}[\cite{MOT}]
\label{thm:MOT}
The Heegaard Floer homologies of three-manifolds and the Ozsv\'ath-Szab\'o four-manifold invariants (mod $2$) are algorithmically computable. 
\end{theorem}

The mod $2$ restriction is a technical point, due to the fact that the link surgery formula has only been proved with $\Z/2$ (rather than with $\Z$) coefficients. The same formula is expected to hold over $\Z$, and in any case the four-manifold invariants mod $2$ suffice to detect many exotic pairs.

\begin{remark}
The first combinatorial description of Heegaard Floer homology was given by Sarkar and Wang in \cite{SarkarWang}, for the hat version $\HFhat$. Theorem~\ref{thm:MOT} extended this to all versions.
\end{remark}

\subsection{Applications.} 
Theorem~\ref{thm:MOT} is a conceptual result, of little use in practice. Interesting closed four-manifolds are represented by complicated links (see the Kirby diagram for $K3$ in Figure~\ref{fig:kirby2}(c)), and therefore need large grid diagrams. If $n$ is the size of the grid diagram, the number of generators of the associated knot Floer complex is $n!$, super-exponential in $n$, which makes computations extremely difficult.

Nevertheless, there are other methods for computing Heegaard Floer invariants. For Heegaard Floer homology (at least in its hat version), an efficient general method is through a theory called {\em bordered Floer homology}. This applies the cut-and-paste principle to three-manifolds by decomposing them into simpler three-manifolds with boundary \cite{LOT, LOTHF}. A similar bordered theory exists for knots \cite{OSbordered, OSbordered2}, which helps to compute knot Floer homology by decomposing the knot into tangles. See \cite{SzaboHFK} for a computer interpretation.

For the moment, bordered theory has not been sufficiently developed to tackle four-manifold invariants directly. The methods that are most effective for that in practice are more {\em ad hoc}: For particular four-manifolds, one can decompose them into handles and take advantage of the properties of the links that we encounter in the process (when we attach $2$-handles, as in a Kirby diagram). For example, in \cite[Section 4]{OSSymp}, Ozsv\'ath and Szab\'o performed  a handle-by-handle calculation of the invariant of the $K3$ surface.

This method was refined in recent work of Levine, Lidman, and Piccirillo \cite{LLP}. There, they construct a new exotic $\CP^2 \# 9\bCP$ handle by handle, and showed it is exotic using the Ozsv\'ath-Szab\'o invariant. The computation involves the relation between knot Floer homology and cobordism maps induced by surgeries, and also makes use of bordered Floer homology. One interesting property of the resulting four-manifold is that it has a free involution. Its quotient by that involution has $\pi_1=\Z/2$ and negative definite intersection form. Here is a corollary.

\begin{theorem}(Levine-Lidman-Piccirillo \cite{LLP}) 
There exists an exotic pair of closed orientable four-manifolds with definite intersection form.
\end{theorem}

This is the first example of this kind. It does not quite answer Question~\ref{q:definite}, because the manifolds are not simply connected (have $\pi_1=\Z/2$). It still shows the power of Heegaard Floer theory, and of the idea of computing four-manifold invariants by using handle decompositions.

In later work, Lidman and Piccirillo used similar methods to construct an exotic $\CP^2 \# 5 \bCP$ \cite{LP}. 

In a different direction, knot Floer homology was used by Juh\'asz and Zemke \cite{JuhaszZemke} to compute the Ozsv\'ath-Szab\'o invariants of four-manifolds obtained by a construction called concordance surgery. (It had not been known how to do the same computation in the Seiberg-Witten setting.)

\section{Skein lasagna modules}
\label{sec:lasagna}

\subsection{Khovanov homology.}
\label{sec:Kh}
In Section~\ref{sec:knots} we briefly mentioned a knot homology theory called Khovanov homology, whose Euler characteristic is the Jones polynomial. We now sketch  its definition.

Recall that the Jones polynomial can be obtained from the Kauffman bracket by a renormalization \eqref{eq:JonesKB}, and that the Kauffman bracket satisfies Equation~\eqref{eq:KB}. By iterating this equation over all $n$ crossings in a diagram $D$, we get a sum of $2^n$ terms, one for each complete resolution of $D$ (where all the crossings are resolved in one of the two ways indicated in Figure~\ref{fig:unorientedskein}). Each complete resolution is an unlink of some number $m$ of components, and therefore its bracket is $(-A^2 - A^{-2})^{m-1}$.

\begin{figure} \centering
{
\fontsize{9pt}{11pt}\selectfont
   \def\svgwidth{3.8in} 
   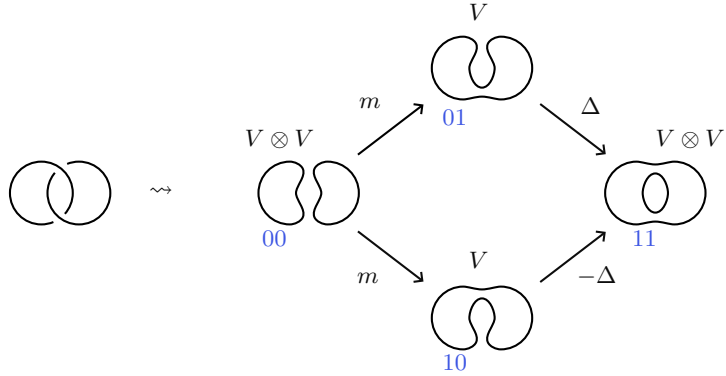
}
\caption{A cube of resolutions and the Khovanov chain complex for a diagram of the Hopf link}
\label{fig:khovanov}
\end{figure}

The construction of Khovanov homology is based on the same idea, with the expression $-A^2 - A^{-2}$ replaced by a free abelian group of rank two, $$V= \text{Span} \ \{1, x\}.$$  We form a hypercube of complete resolutions as in Figure~\ref{fig:khovanov}, where the vertices are labeled by $\eps=(\eps_1, \dots, \eps_n) \in \{0,1\}^n$. Then, we form a chain complex where at each vertex we place a tensor product $V^{\otimes m}$. We take the direct sum of all the groups associated to a value of $|\eps|=\eps_1 + \dots +\eps_n$. These are our chain groups. 

The differential $d$ relates the groups at $\eps$ to those at $\eps'$ where $\eps'$ differs from $\eps$ at a single crossing, say $j \in \{1, \dots, n\}$, so that $\eps_j=0$ but $\eps'_j =1$. The resolution at $\eps'$ is obtained from that at $\eps$ by either joining two circles or splitting a circle into two. 

If two circles get joined together, we define the contribution of $d$ using the multiplication map
$$ m: V \otimes V \to V$$
$$1^2 =1, \ 1x = x1 = x, \ x^2 =0,$$
tensored with the identity on all the other factors (for unchanged circles). 

If one circle gets split in two, we use the comultiplication
$$\Delta: V \to V \otimes V$$
$$ \Delta(1)= 1 \otimes x + x \otimes 1, \ \ \Delta(x) = x \otimes x,$$
and again tensor with the identity on the remaining factors. 

We also multiply each term in the differential by a sign $(-1)^{\eps_1 + \dots + \eps_{j-1}}$. One can check that the resulting map $d$ satisfies $d^2=0$, and thus defines a chain complex. The homology of that complex is Khovanov homology, $\Kh(K)$. (We can further equip it with two gradings as in \eqref{eq:Kh}, but we will not discuss those here.)

The definition above is purely algebraic, and Khovanov homology is deeply connected to modern representation theory. However, it also shares some properties with the analytically defined Floer homologies from Subsection~\ref{sec:3mi}. Specifically, it is functorial under link cobordisms, as follows. Suppose we have a smoothly embedded surface $\Sigma \subset [0,1]\times \R^3$ with boundary $\del \Sigma = (-L_0) \sqcup L_1$, for some links $L_0 \subset \{0\} \times \R^3$ and $L_1 \subset \{1\} \times \R^3$. Then, we have a well-defined map
$$ \Kh(\Sigma): \Kh(L_0) \to \Kh(L_1),$$
and composition of cobordisms corresponds to composition of maps. See \cite{Jacobsson, KhTangleCob, BN,  Blanchet}. 

This TQFT property lies at the heart of the construction of four-manifold invariants from Khovanov homology that will be presented in Subsection~\ref{sec:slm}.

\subsection{The skein module.}
The four-manifold invariants were inspired by a skein-theoretic construction one dimension lower, which we now review. Given a three-manifold $Y$, the {\em Kauffman bracket skein module} of $Y$, denoted $\KBSM(Y)$, is generated over $\Z[A, A^{-1}]$ by all framed links in $Y$ modulo the local relations \eqref{eq:KB1}, \eqref{eq:KB} from the definition of the Kauffman bracket; see \cite{Przytycki, Turaev}. Each framed link has an invariant, its class $\langle \Link \rangle \in \KBSM(Y)$.

In the case of $Y=S^3$, one can easily check that $\KBSM(S^3) = \Z[A, A^{-1}]$, with a generator being the empty link. The class of a framed link $\Link \subset S^3$ is simply its Kauffman bracket. 

For general three-manifolds, one can view $\langle \Link \rangle$ as a generalization of the Jones polynomial (up to a normalization factor). However, while the Kauffman bracket skein module is easy to define, it is hard to compute. There is an extensive literature devoted to its computation. It took many years to even establish a basic finite dimensionality property: in \cite{GJS}, Gunningham, Jordan and Safronov proved a conjecture of Witten, that $\KBSM(Y) \otimes \C(A)$ is finite dimensional over $\C(A)$. It is expected that $\KBSM(Y) \otimes \C(A)$ is related to the $\mathit{SL}(2, \C)$ Floer homology of three-manifolds defined in \cite{AM}.

\subsection{A new invariant of four-manifolds.}
\label{sec:slm}
We now describe an invariant introduced by Morrison, Walker and Wedrich in \cite{MWW}. It is associated to a link in the boundary of any four-manifold; by taking the link to be empty, it also gives an invariant of the four-manifold. The invariant generalizes Khovanov homology in the same sense that the link invariants in the Kauffman bracket skein module generalize the Jones polynomial. The idea is that, instead of imposing the local skein relations in 3D, we impose local relations from cobordism maps in 4D.

In the set-up of \cite{MWW}, they work with a differently normalized version of Khovanov homology, that depends on the framing of the link; for simplicity, we still denote it by $\Kh$. This is functorial under framed surface cobordisms in $[0,1] \times \R^3$.  

In fact, we need a more refined cobordism property, proved in \cite{MWW}. This says that whenever we have a collection of four-balls $B_1, \dots, B_k$ in the interior of another four-ball $B$, and we are given a framed surface $$\Sigma \subset B \setminus (\cup_{i=1}^k B_i)$$
with boundaries $\Link_i \subset \del B_i$ and $\Link \subset \del B$, there is a well-defined cobordism map
$$ \Kh(\Sigma): \bigotimes_{i=1}^k \Kh(\Link_i) \to \Kh(\Link).$$
To define it, we remove a basepoint $z_i$ from each $\del B_i$, a basepoint $z$ from $\del B$, as well as an embedded tree connecting $z_i$ and $z$ and disjoint from $\Sigma$. The result is that $\Sigma$ now lives in $[0,1] \times \R^3$, so we can use the ordinary Khovanov map. The extra input from \cite{MWW} is that this map does not depend on what basepoints and tree we removed.  

Now let $X$ be a compact four-manifold with boundary $\del X=Y$. (Recall from Corollary~\ref{cor:3bdry} that every $Y$ appears this way.) Let $\Link$ be a framed link in $Y$. 

\begin{figure} \centering
{
 \includegraphics[scale=.32]{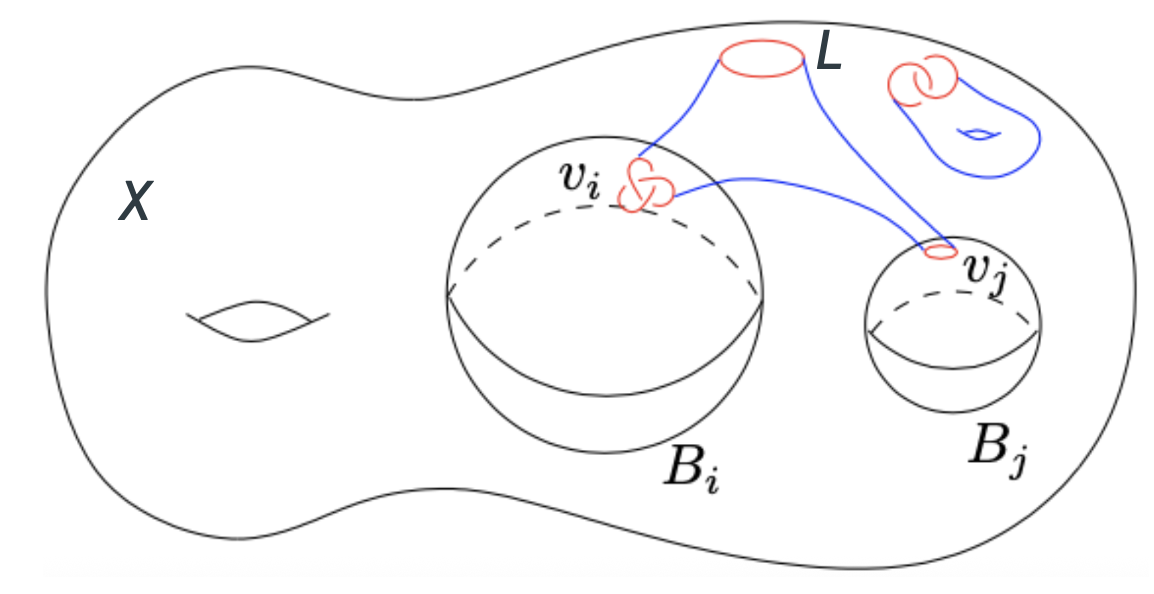}
}
\caption{A lasagna filling}
\label{fig:lasagna1}
\end{figure}

\begin{definition}
A {\em lasagna
filling} $F=(\Sigma, \{(B_i,\Link_i,v_i)\})$ of $X$ with boundary $\Link$ consists of
\begin{itemize}
\item A finite collection of disjoint $4$-balls $B_i$ (called {\em input balls})
embedded in the interior of $X$;
\item A framed surface $\Sigma$ properly embedded in $X \setminus
\cup_i B_i$, meeting $\del X$ in $\Link$ and meeting each $\del B_i$ in a link
$\Link_i$; and
\item for each $i$, an element $v_i \in \Kh( \Link_i).$
\end{itemize}
See Figure~\ref{fig:lasagna1}.
\end{definition}

\begin{figure} \centering
{
\fontsize{9pt}{11pt}\selectfont
   \def\svgwidth{5in} 
   \resizebox{0.8\textwidth}{!}{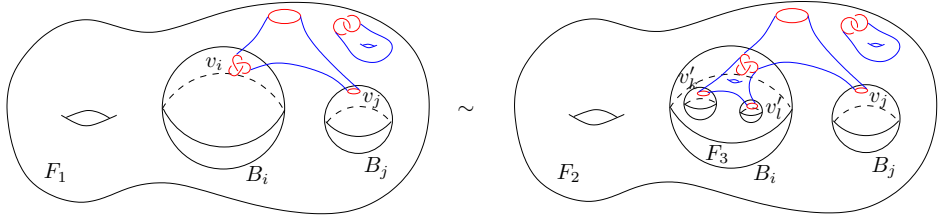}
}
\caption{Equivalence of lasagna fillings}
\label{fig:lasagna2}
\end{figure}

\begin{definition}
The {\em skein lasagna module} of $\Link$ relative to $X$ is
$$ \Sz(X; L) := \Z\{ \text{lasagna fillings $F$ of $X$ with boundary
$\Link$}\}/\sim$$ where $\sim$ is the transitive and linear closure of the following
relations: 
\begin{itemize}
\item Linear combinations of lasagna fillings are set to be multilinear in the
elements $v_i$;
\item  $F_1$ and $F_2$ are set to be equivalent
if  $F_1$ has an input ball $B_i$ with label $v_i$, and $F_2$ is obtained from
$F_1$ by replacing $B_i$ with another lasagna filling $F_3=(\Sigma', \{(B'_i,\Link'_i,v'_i)\})$ of a four-ball such
that $v_i=\Kh(\Sigma')(\otimes v'_i)$, followed by an isotopy rel $\del X$. See Figure~\ref{fig:lasagna2}.
\end{itemize}
In particular, if $X$ is any compact four-manifold (with or without boundary), we let $\Sz(X) := \Sz(X; \emptyset)$.
\end{definition}

The name {\em skein lasagna module} is inspired by some earlier terminology due to Jones: One dimension lower, one considers ``spaghetti and meatballs" pictures consisting of linked paths joining three-balls inside another three-ball. In dimension four, the lasagna fillings are two-dimensional, hence their name. 

The procedure we described above upgrades a link invariant to a four-manifold invariant, and is quite general. It can be applied not only to Khovanov homology but to any other link homology satisfying similar properties. In  \cite{MWW} the authors apply it to the Khovanov-Rozansky homologies associated to the Lie algebras $\mathfrak{gl}_N$ for any $N$. (The case $N=2$ is Khovanov homology.) Furthermore, in  \cite{ChenFloer}, Chen upgrades knot Floer homology in an analogous way, and defines an invariant called the {\em Floer lasagna module}.

Going back to the skein lasagna modules from Khovanov homology, let us discuss a few of their basic properties. 
First, they are (by construction) functorial under cobordisms, in the following sense. Given a four-manifold $W$ with boundary $Y$, a cobordism $Z$ from $Y$ to $Y'$, and a properly embedded framed surface $\Sigma \subset Z$ with boundaries $\Link \subset Y$ and $\Link' \subset Y'$, there is a natural map
$$ \Psi_{Z, \Sigma}: \Sz(W; \Link) \to \Sz(W \cup Z; \Link')$$
given by attaching $\Sigma$ to lasagna fillings.

Second, the bigrading on Khovanov homology descends to one on skein lasagna modules, and there is a further decomposition according to the relative homology class of the lasagna filling:
$$ \Sz(W; \Link) =\bigoplus_{i, j \in \Z} \bigoplus_{\alpha \in \del^{-1}([L])} \Sz_{i, j}(W; \Link, \alpha).$$
Here, $\Link=(L, \lambda)$ and $\del^{-1}([L]) \subset H_2(W; L)$ is the preimage of the fundamental class $[L] \in H_1(\Link)$ under the boundary homomorphism. A lasagna filling can be further filled with surfaces inside the input ball to produce the class $\alpha$.

\subsection{Computational methods.}
As in the case of the 3D skein modules, the 4D skein lasagna modules are (relatively) easy to define, but hard to compute. 

The simplest calculation can be done in the case $X=D^4$. Then $\Link \subset \del D^4=S^3$ has a Khovanov homology itself, $\Kh(\Link)$, and there is a natural map 
$$\Sz(X; \Link) \to \Kh(\Link), \ \ [(\Sigma, \{(B_i,\Link_i,v_i)\})] \to \Kh(\Sigma)(\otimes_i  v_i).$$
It is not hard to see that this is an isomorphism, so $\Sz(X; \Link) \cong \Kh(\Link)$. Thus, skein lasagna modules   generalize Khovanov homology.

To do other calculations, it is helpful to decompose the four-manifold into handles as in Section~\ref{sec:kirby}. Surgery formulas have been developed in \cite{MN} and \cite{ManWW},  describing how the skein lasagna module behaves under a handle attachment. The formulas  for $k$-handles are getting more complicated as $k$ gets smaller.

To begin, attaching a $4$-handle does not change $\Sz(W; \Link)$. For example, we have $\Sz(D^4) = \Sz(S^4) =\Z$.

Next, let us consider a $3$-handle attachment to $(W; \Link)$, along a sphere $S^2 \subset Y=\del W$ (disjoint from $\Link$). Let the result be $W'$. Break the sphere into two hemispheres $\Delta_+$ and $\Delta_-$, and push their interiors slightly inside a cylinder $I \times Y$, where $I=[0,1]$. Consider the difference of two cobordism maps
$$ f=\Psi_{I \times Y, \Delta_+} - \Psi_{I \times Y, \Delta_-}: \Sz(W; \Link) \to \Sz(W; \Link),$$
where we identified $W$ with $W \cup (I \times Y)$. 

\begin{theorem}[Theorem 3.7 in \cite{ManWW}]
If $W'$ is obtained from $W$ by attaching a $3$-handle as above, then
$$ \Sz(W'; \Link)\cong \Sz(W; L)/\operatorname{im}(f).$$
\end{theorem}

We move on to the $2$-handle formula. Suppose $W'$ is obtained from $W$ by attaching $2$-handles along a framed link $\K \subset Y=\del W$. We are also given a framed link $\Link \subset Y$ disjoint from $\K$ (but possibly linked with it). We let $\K(r_1, r_2)$ denote the cable of $\K$ consisting of $r_1+r_2$ parallel push-offs of $\K$ (according to its framing), where $r_1$ of them are oriented the same way as $\K$ and $r_2$ are oriented in the opposite direction. Let also $U$ be the unknot. For each $r_1$ and $r_2$, there is a cobordism $Z$ from the disjoint union $\K(r_1, r_2) \cup L \sqcup U$ to $\K(r_1+1, r_2 +1) \cup L$, given by removing a small disk from an annulus that goes between two oppositely oriented strands of $\K(r_1+1, r_2 +1)$.

\begin{theorem}[Theorem 1.1 in \cite{MN}; Theorem 3.2 in \cite{ManWW}] 
\label{thm:MN}
If $W'$ is obtained from $W$ by attaching a $2$-handle as above, then
\begin{equation}
\label{eq:s2h}
 \Sz(W'; \Link) \cong \bigoplus_{r_1, r_2 \in \N} \Sz(W; \K(r_1, r_2) \cup \Link)/\sim, 
 \end{equation}
where the equivalence $\sim$ is generated by the following:
\begin{itemize}
\item permuting the strands of $\K(r_1, r_2)$ in an orientation-preserving way,
\item $\Psi_Z(v \otimes 1) \sim 0$,
\item $\Psi_Z(v \otimes X) \sim v$, 
\end{itemize}
for all $v \in \Sz(W; \K(r_1, r_2) \cup \Link).$ (The skein lasagna module gets tensored with $\Kh(U)$ under the split disjoint union with $U$, and the elements $1$ and $X$ are the generators of $\Kh(U)$.)
\end{theorem}

The first version of Theorem~\ref{thm:MN} was proved in \cite{MN}, for the case where $W=D^4$, so that $W'$ is a $2$-handlebody (a four-manifold made by attaching $2$-handles to $D^4$). Suppose also that $\Link =0$. In that situation the right hand side of Equation~\eqref{eq:s2h} involves the Khovanov homologies of the cables $\K(r_1, r_2)$. To compute $\Sz(W')$, one needs to understand all these homologies, which is difficult in general. Nevertheless, it can be done for example when $\K$ is the $0$-framed unknot, so that its cables are unlinks. This led to a calculation of $\Sz(S^2 \times D^2)$.

By taking $\K$ to be the unknot with $\pm 1$ framings, and then using the fact that a four-handle does not change the skein lasagna module, in \cite{MN} we obtained partial computations of $\Sz(\CP)$ and $\Sz(\bCP)$. This sufficed to prove that these two groups are different, so the skein lasagna module is sensitive to orientation.

The methods in \cite{MN} were further refined in work of Sullivan-Zhang \cite{SullivanZhang} and Ren-Willis \cite{RenWillis}. For example, they proved that
$$\Sz(\CP^2)=\Sz(S^2 \times S^2) = 0.$$

There is also a formula for attaching $1$-handles; see \cite[Theorem 4.7]{ManWW}. We will not state it here, as it is quite complicated. It suffices to say that in the case of a single $1$-handle, and a framed link $\Link$ in the boundary of the resulting $S^1 \times D^3$, the formula involves the Hochschild homology of a category of tangles for a bimodule associated to $\Link$. One application is the following calculation. For the link $\Link_p \subset S^1 \times S^2$ that consists of $2p$ parallel longitudes $S^1 \times \{\text{pt}\}$, half of which are oriented one way and half the other way:
$$ \Sz(S^1 \times D^3, \Link_p) \otimes \Q \cong \begin{cases}
\Q & \text{if } p=0,\\ 
\Q^4 & \text{if } p=1,\\
\Q^\infty & \text{if } p=\infty. 
\end{cases}$$

Altogether, the formulas in \cite{MN} and \cite{ManWW} give a general description of $\Sz(X; L)$ from handle decompositions. These look quite different from the corresponding formulas for Heegaard Floer homology in Section~\ref{sec:surgeryHF}, but have the same effect of enabling the computation of four-manifold invariants in terms of link invariants.

\subsection{Detection of exotic pairs.} 
\begin{figure} \centering
{
\fontsize{11pt}{11pt}\selectfont
   \def\svgwidth{3.2in} 
   %% Creator: Inkscape 1.3.2 (091e20e, 2023-11-25), www.inkscape.org
%% Accompanies image file 'exotic.pdf' (pdf, eps, ps)
\begingroup%
  \makeatletter%
  \providecommand\color[2][]{%
    \errmessage{(Inkscape) Color is used for the text in Inkscape, but the package 'color.sty' is not loaded}%
    \renewcommand\color[2][]{}%
  }%
  \providecommand\transparent[1]{%
    \errmessage{(Inkscape) Transparency is used (non-zero) for the text in Inkscape, but the package 'transparent.sty' is not loaded}%
    \renewcommand\transparent[1]{}%
  }%
  \providecommand\rotatebox[2]{#2}%
  \newcommand*\fsize{\dimexpr\f@size pt\relax}%
  \newcommand*\lineheight[1]{\fontsize{\fsize}{#1\fsize}\selectfont}%
  \ifx\svgwidth\undefined%
    \setlength{\unitlength}{377.78197762bp}%
    \ifx\svgscale\undefined%
      \relax%
    \else%
      \setlength{\unitlength}{\unitlength * \real{\svgscale}}%
    \fi%
  \else%
    \setlength{\unitlength}{\svgwidth}%
  \fi%
  \global\let\svgwidth\undefined%
  \global\let\svgscale\undefined%
  \makeatother%
  \begin{picture}(1,0.46153559)%
    \lineheight{1}%
    \setlength\tabcolsep{0pt}%
    \put(0,0){\includegraphics[width=\unitlength,page=1]{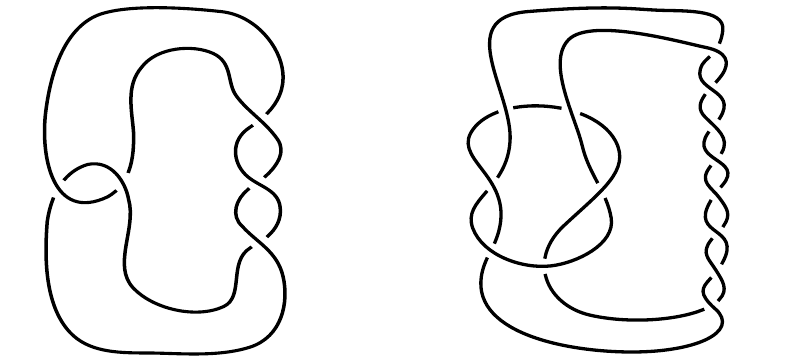}}%
    \put(0.01953602,0.00888119){\color[rgb]{0,0,0}\makebox(0,0)[lt]{\lineheight{1.25}\smash{\begin{tabular}[t]{l}$K_1$\end{tabular}}}}%
    \put(0.59218776,0.01116567){\color[rgb]{0,0,0}\makebox(0,0)[lt]{\lineheight{1.25}\smash{\begin{tabular}[t]{l}$K_2$\end{tabular}}}}%
    \put(0.35149049,0.43375206){\color[rgb]{0,0,0}\makebox(0,0)[lt]{\lineheight{1.25}\smash{\begin{tabular}[t]{l}$-1$\end{tabular}}}}%
    \put(0.93543336,0.43146759){\color[rgb]{0,0,0}\makebox(0,0)[lt]{\lineheight{1.25}\smash{\begin{tabular}[t]{l}$-1$\end{tabular}}}}%
  \end{picture}%
\endgroup%

}
\caption{An exotic pair of traces}
\label{fig:exotic}
\end{figure}
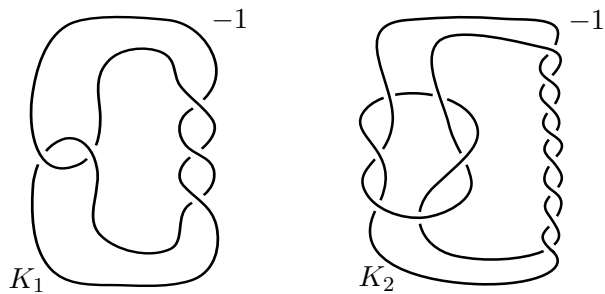

A striking application of skein lasagna modules is a proof, by Ren and Willis, that these invariants can distinguish exotic smooth structures on some compact four-manifolds with boundary \cite{RenWillis}. The example they consider is a pair of traces of $-1$ surgeries on knots: One knot is $K_1=-5_2$ and the other is the pretzel knot $K_2=P(3, -3, -8)$; see Figure~\ref{fig:exotic}. The traces $W_1 =X(K_1, -1)$ and $W_2=X(K_2, -1)$ have the same boundary $S^3(K_1, -1) = S^3(K_2, -1)$, are simply connected, and have the same homology and intersection form. This implies that they are homeomorphic by \cite{Boyer}. They are not diffeomorphic because their skein lasagna modules are different. For example,
$$ \Sz_{0, q} (W_1; 1) \otimes \Q \cong \begin{cases}
\Q & \text{if $q=1,3$},\\
0 & \text{otherwise,}
\end{cases}$$
whereas
$$ \Sz_{0, q} (W_2; 1) \otimes \Q \supseteq \begin{cases}
\Q & \text{if $q=-1,1$},\\
0 & \text{otherwise.}
\end{cases}$$
These calculations are based on the $2$-handle formula in \cite{MN}, and on a partial understanding of Khovanov homology of the cables of $K_1$ and $K_2$, using special properties of these knots. In the case of $K_1$, the key property is that it admits a diagram with only positive crossings. In the case of $K_2$, a useful property is that it is slice. (See Subsection~\ref{sec:slice} below for a definition of sliceness.)

It had been known before that $W_1$ and $W_2$ are an exotic pair, by the work of Akbulut \cite{Akbulut} that used gauge theory. Ren and Willis gave the first analysis-free proof. In their paper they also exhibited some new exotic examples: detectable using skein lasagna modules, but not (for the moment) using gauge theory or Heegaard Floer homology. All their examples have non-empty boundary, and the following is still an open problem.

\begin{question}
Can skein lasagna modules detect exotic smooth structures on some closed four-manifolds?
\end{question}

In a related direction, Sullivan \cite{Sullivan} showed that a variant of skein lasagna modules detects some exotically knotted pairs of surface embeddings in the four-ball. 

Recently, Nahm \cite{Nahm} gave a more direct way of detecting some exotic compact four-manifolds with boundary, using the cobordism maps on Khovanov homology and without needing skein lasagna modules. The simplest of his examples are the complements of the exotic disks in $D^4$ studied in \cite{AkbulutZeeman, HaydenSundberg}.

\section{Probing four-manifolds with knots}
\label{sec:probing}

\subsection{Sliceness.} \label{sec:slice} A knot $K \subset S^3$ is the unknot if and only if it bounds a smoothly embedded disk. By contrast, if we go to four dimensions by including $K \subset S^3 \subset S^4$, then every knot $K$ bounds a smoothly embedded disk in $S^4$: we can change the crossings at will and isotope the knot into the unknot in $S^4$, extend this isotopy to an ambient isotopy of $S^4$ itself, and follow the disk bounded by the unknot in reverse under the ambient isotopy. 

There is, however, an intermediate problem that is more interesting. Instead of allowing the disk to run freely into $S^4$, we ask for it to be properly embedded in one hemisphere $D^4$, whose boundary $S^3$ contains the knot. 

\begin{definition}
A knot $K \subset S^3 = \del D^4$ is called {\em (smoothly) slice} if it bounds a smoothly, properly embedded disk in $D^4$.
\end{definition}

\begin{figure} \centering
{
 \includegraphics[scale=.37]{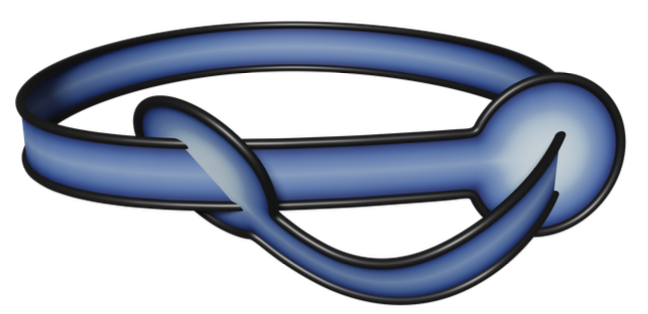}
}
\caption{A slice knot}
\label{fig:slice}
\end{figure}

\begin{figure} \centering
{
\fontsize{9pt}{11pt}\selectfont
   \def\svgwidth{4.5in} 
   %% Creator: Inkscape 1.3.2 (091e20e, 2023-11-25), www.inkscape.org
%% Accompanies image file 'Ribbon.pdf' (pdf, eps, ps)
\begingroup%
  \makeatletter%
  \providecommand\color[2][]{%
    \errmessage{(Inkscape) Color is used for the text in Inkscape, but the package 'color.sty' is not loaded}%
    \renewcommand\color[2][]{}%
  }%
  \providecommand\transparent[1]{%
    \errmessage{(Inkscape) Transparency is used (non-zero) for the text in Inkscape, but the package 'transparent.sty' is not loaded}%
    \renewcommand\transparent[1]{}%
  }%
  \providecommand\rotatebox[2]{#2}%
  \newcommand*\fsize{\dimexpr\f@size pt\relax}%
  \newcommand*\lineheight[1]{\fontsize{\fsize}{#1\fsize}\selectfont}%
  \ifx\svgwidth\undefined%
    \setlength{\unitlength}{433.62029308bp}%
    \ifx\svgscale\undefined%
      \relax%
    \else%
      \setlength{\unitlength}{\unitlength * \real{\svgscale}}%
    \fi%
  \else%
    \setlength{\unitlength}{\svgwidth}%
  \fi%
  \global\let\svgwidth\undefined%
  \global\let\svgscale\undefined%
  \makeatother%
  \begin{picture}(1,0.13489623)%
    \lineheight{1}%
    \setlength\tabcolsep{0pt}%
    \put(0,0){\includegraphics[width=\unitlength,page=1]{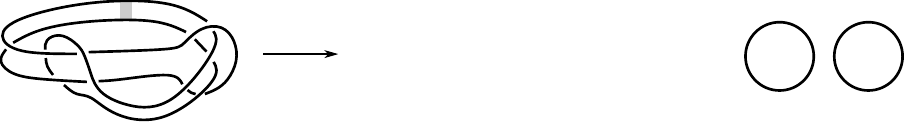}}%
    \put(0.30603264,0.0900938){\makebox(0,0)[lt]{\lineheight{1.25}\smash{\begin{tabular}[t]{l}band\end{tabular}}}}%
    \put(0.74014413,0.08836552){\makebox(0,0)[lt]{\lineheight{1.25}\smash{\begin{tabular}[t]{l}$\cong$\end{tabular}}}}%
    \put(0,0){\includegraphics[width=\unitlength,page=2]{Ribbon.pdf}}%
  \end{picture}%
\endgroup%

}
\caption{Attaching a band}
\label{fig:band}
\end{figure}

An example of a non-trivial slice knot is the one in Figure~\ref{fig:slice}. It bounds an immersed disk in $\R^3$ as in the picture. We can think of the intensity of the color as the fourth dimension, and thus the disk becomes embedded in $D^4$. More generally, when an immersed disk intersects itself only along intervals as in the local models seen in  Figure~\ref{fig:slice}, we say it is a ribbon disk, and its boundary knot is called a {\em ribbon knot}. We can get rid of each ribbon singularity by attaching a band to the knot as in Figure~\ref{fig:band}. One can easily see that a knot is ribbon if and only if there exist some number of $k$ bands such that attaching them produces the $(k+1)$-component unlink.

Observe that every ribbon knot is slice. The converse is a famous open problem.

\begin{conjecture}[The slice-ribbon conjecture \cite{Fox}]
Every slice knot is ribbon.
\end{conjecture}
 
 Sliceness is a four-dimensional property of knots.  Many problems in three-dimensional topology are decidable: for example, there exist algorithms for deciding whether two knots are isotopic, or whether two three-manifolds are homeomorphic; see for example \cite{Hemion, ScottShort}. In four-dimensional topology, on the other hand, many similar problems (such as the homeomorphism problem for four-manifolds) are undecidable \cite{Markov}. When it comes to sliceness or ribbon-ness (which lie at the interface of three- and four-dimensional topology), we do not even know if these problems are decidable. 

While there is no known algorithm for detecting slice knots, in practice one can go quite far using two kinds of tools:
\begin{enumerate}
\item[(a)] {\em obstructions} for a knot to be slice; or, equivalently, necessary conditions on knot invariants for this property to hold. For example, we have the Fox-Milnor condition, that the Alexander polynomial should be of the form $f(q)f(q^{-1})$ for some Laurent polynomial $f(q)$ \cite{FoxMilnor}. Other knot invariants are known to vanish for slice knots: Examples include topological invariants (such as the signature $\sigma$), some coming from gauge theory or Heegaard Floer homology (having names such as $\tau$, $\epsilon$, $\nu$, $\delta$, etc.), and some from Khovanov homology and related algebraic theories (the best known example of these is an invariant denoted $s$ that was constructed by Rasmussen \cite{RasmussenSlice});

\item[(b)] {\em constructions} of slice disks. In most cases, one simply looks for a ribbon disk in $\R^3$, and one can do so by searching for bands that transform the knot into the unlink. There are computer programs that do this, some involving machine learning (Bayesian optimization) \cite{GHMR}, and some classical but still very effective \cite{DunfieldGong}. 
\end{enumerate}

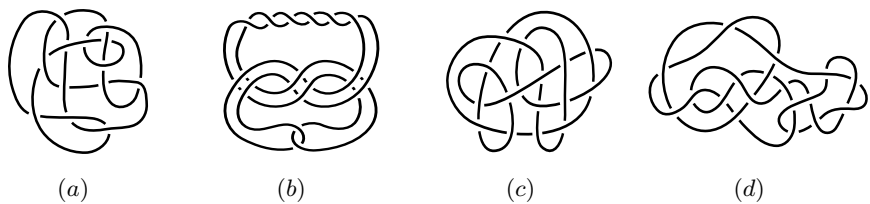
\begin{figure} \centering
{
\fontsize{9pt}{11pt}\selectfont
   \def\svgwidth{4.5in} 
   %% Creator: Inkscape 1.3.2 (091e20e, 2023-11-25), www.inkscape.org
%% Accompanies image file 'unknown.pdf' (pdf, eps, ps)
\begingroup%
  \makeatletter%
  \providecommand\color[2][]{%
    \errmessage{(Inkscape) Color is used for the text in Inkscape, but the package 'color.sty' is not loaded}%
    \renewcommand\color[2][]{}%
  }%
  \providecommand\transparent[1]{%
    \errmessage{(Inkscape) Transparency is used (non-zero) for the text in Inkscape, but the package 'transparent.sty' is not loaded}%
    \renewcommand\transparent[1]{}%
  }%
  \providecommand\rotatebox[2]{#2}%
  \newcommand*\fsize{\dimexpr\f@size pt\relax}%
  \newcommand*\lineheight[1]{\fontsize{\fsize}{#1\fsize}\selectfont}%
  \ifx\svgwidth\undefined%
    \setlength{\unitlength}{469.20336626bp}%
    \ifx\svgscale\undefined%
      \relax%
    \else%
      \setlength{\unitlength}{\unitlength * \real{\svgscale}}%
    \fi%
  \else%
    \setlength{\unitlength}{\svgwidth}%
  \fi%
  \global\let\svgwidth\undefined%
  \global\let\svgscale\undefined%
  \makeatother%
  \begin{picture}(1,0.2233023)%
    \lineheight{1}%
    \setlength\tabcolsep{0pt}%
    \put(0,0){\includegraphics[width=\unitlength,page=1]{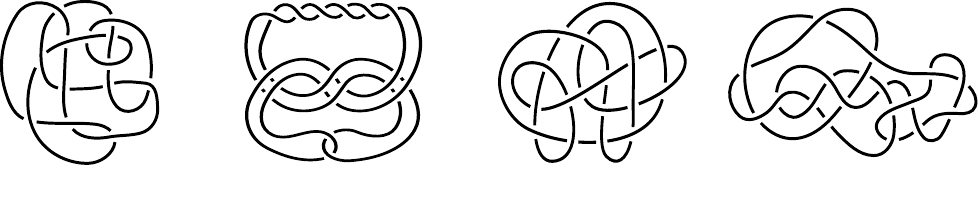}}%
    \put(0.05721671,0.00560237){\color[rgb]{0,0,0}\makebox(0,0)[lt]{\lineheight{1.25}\smash{\begin{tabular}[t]{l}$(a)$\end{tabular}}}}%
    \put(0.31236554,0.00560237){\color[rgb]{0,0,0}\makebox(0,0)[lt]{\lineheight{1.25}\smash{\begin{tabular}[t]{l}$(b)$\end{tabular}}}}%
    \put(0.57905418,0.00560237){\color[rgb]{0,0,0}\makebox(0,0)[lt]{\lineheight{1.25}\smash{\begin{tabular}[t]{l}$(c)$\end{tabular}}}}%
    \put(0.84388617,0.00560237){\color[rgb]{0,0,0}\makebox(0,0)[lt]{\lineheight{1.25}\smash{\begin{tabular}[t]{l}$(d)$\end{tabular}}}}%
  \end{picture}%
\endgroup%

}
\caption{Some knots whose sliceness is unknown: (a) $13n65$; (b) the Whitehead double of the left-handed trefoil; (c) $16n68728$; (d) $18nh_{00000601}$}
\label{fig:unknown}
\end{figure}

By combining obstructive and constructive methods, Dunfield and Gong \cite{DunfieldGong} went through the list of all $\approx \! 350$ million prime knots with up to $19$ crossings, and showed that approximately $99.5\%$ of them are not slice, and $0.5\%$ are slice. This left only a small percentage ($0.003\%$, i.e. $\approx \! 11,400$ knots) for which they could not determine sliceness. The smallest knots of unknown sliceness have $13$ crossings; an example is in Figure~\ref{fig:unknown} (a). 

One famous knot of unknown sliceness is the Whitehead double of the left-handed trefoil, which  is shown in Figure~\ref{fig:unknown} (b). Two other notable knots, whose sliceness had been not known for a long time, had recently been proven to not be slice: the Conway knot \cite{Pic2} and the $(2,1)$-cable of the figure-eight \cite{DKMPS}.

\subsection{Sliceness and H-sliceness in four-manifolds.}
There are versions of the slice property relative to any closed four-manifold $X$. Let us denote by $X^\circ$ the complement of an open four-ball in $X$. This is a compact four-manifold with boundary $S^3$.

\begin{definition}
 $(a)$ A knot $K \subset S^3=\del X^\circ$ is called {\em slice in $X$} if it bounds a smoothly, properly embedded disk $\Delta \subset X^\circ$. Such a disk has a relative homology class  $[\Delta] \in H_2(X^\circ, \del X^\circ) \cong H_2(X)$. 
 
 $(b)$ If a disk $\Delta$ exists such that $[\Delta]=0$, the knot $K$ is called {\em H-slice in $X$}. 
 
 $(c)$ If a disk $\Delta$ exists such that $[\Delta]^2 = -k \in \Z$, the knot $K$ is called {\em $k$-slice in $X$}.
\end{definition}

When $X=S^4$, we recover the usual notion of slice, which in this case is equivalent to being H-slice, or $0$-slice.

There are again many obstructions and constructions that can be used to study the different sliceness properties in four-manifolds. We refer to \cite[Section 2]{MMP} for some examples. For now let us mention that the left handed and the right handed trefoils, while not slice in $S^4$, are slice in $\CP^2$; the right handed trefoil is even H-slice in $\CP^2$. There exist knots that are not slice in $\CP^2$, whereas every knot is slice in $S^2 \times S^2$. It is unknown whether every knot is slice in the $K3$ surface, but this property is known to hold for all knots that can be unknotted by at most $21$ crossing changes \cite{MarengonMihajlovic}. 
 
\subsection{A program to find new exotica.}
Rather than fixing a four-manifold $X$ and asking which knots are slice (or H-slice, or $k$-slice) in it, we can reverse the question: Given a knot $K \subset S^3$, what do these properties tell us about the four-manifold? In particular, can they detect exotic pairs?

In the case of H-sliceness, an application of Seiberg-Witten theory shows that the answer is affirmative:
\begin{theorem}[Corollary 1.4 in \cite{MMP}]
There exist smooth, homeomorphic four-manifolds $X$ and $X'$ and a knot $K \subset S^3$ that is H-slice in $X$ but not in $X'$. As an example, one can take $$X=\#3\CP^2 \# 20 \bCP, \ \ X' = K3 \# \bCP,$$ and $K$ to be the right-handed trefoil. 
\end{theorem}

One can also find such examples for $k$-sliceness. The corresponding problem for sliceness is still open, although some progress in that direction exists: Lidman and Piccirillo \cite{LP} showed that there exist four-manifolds $X$ and $X'$ with the same cohomology ring, and a knot that is slice in $X$ but not in $X'$. (Their manifolds are not homeomorphic though, because they have different fundamental groups.)

\begin{question}
Does there exist an exotic pair of closed four-manifolds $(X, X')$ and a knot $K \subset S^3$ that is slice in $X$ but not in $X'$?
\end{question}

In principle, one could use this idea to try to disprove the smooth four-dimensional Poincar\'e conjecture (SPC4). There exist many potential counterexamples (manifolds known to be homeomorphic to $S^4$, but not known to be diffeomorphic to it; i.e., homotopy $4$-spheres not known to be standard). To run this program, one would also need:
\begin{enumerate}
\item [(a)] an obstruction to sliceness in $S^4$ that does not apply to sliceness in other homotopy $4$-spheres, and 
\item[(b)] a construction of knots that are slice in a homotopy $4$-sphere.
\end{enumerate} 

With regard to (a), most slice obstructions from classical topology, gauge theory, or Heegaard Floer homology apply just as well to sliceness in homotopy $4$-spheres as in $S^4$. More promising are the obstructions from Khovanov homology, in particular Rasmussen's $s$-invariant. This vanishes for knots that are slice in the usual sense (in $S^4$), and it is not known whether the same holds in all other homotopy $4$-spheres. However, a  popular construction of potential counterexamples to SPC4 is by Gluck twists \cite{Gluck}, and it can be proved  that $s(K)=0$ for knots $K$ that are slice in a homotopy $4$-sphere obtained from $S^4$ by a Gluck twist  \cite{MMSW}. Thus, one needs to avoid this construction.

With regard to (b), one idea is to consider homotopy $4$-spheres without $1$-handles, and band together the components of the attaching link for the $2$-handles. An early attempt to disprove SPC4 in this way was in \cite{FGMW}, but the homotopy $4$-spheres in question turned out to be standard \cite{Akbulut4}.

\begin{figure} \centering
{
\fontsize{9pt}{11pt}\selectfont
   \def\svgwidth{3.6in} 
 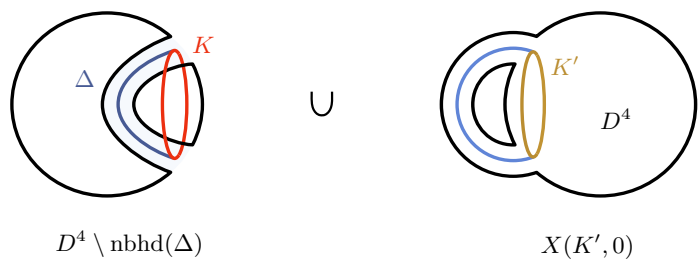
}
\caption{A homotopy $4$-ball obtained by gluing a slice disk complement to the trace of a $0$-surgery}
\label{fig:htpy4}
\end{figure}

A related idea is to consider pairs of knots $(K, K')$ with the same $0$-surgery: $S^3(K, 0) = S^3(K', 0)$; these are called {\em $0$-friends}. Suppose we could find such a pair such that $K$ is slice and $K'$ is not; for example, the sliceness of $K'$ could be obstructed by the non-vanishing of the Rasmussen invariant. Then, one could consider the smooth disk $\Delta \subset D^4$ with boundary $K$, take the complement of a tubular neighborhood of it, and glue it to the trace of the other knot:
$$W = (D^4 \setminus \operatorname{nbhd}(\Delta)) \cup_{S^3(K, 0)} X(K', 0).$$ 
It is easy to check that $W$ is a homotopy $4$-sphere, and the knot $K'$ is slice in it by construction---because removing a ball from $X(K', 0) \subset W$ exhibits a $2$-handle whose core disk has boundary $K'$. (See Figure~\ref{fig:htpy4}.) Since $K'$ is not slice in $S^4$, we would deduce that $W$ is not diffeomorphic to $S^4$, hence disproving SPC4.

As of now, no pairs $(K, K')$ with the above properties have been found. A general construction of $0$-friend pairs $(K, K')$  was given by the author and Piccirillo in \cite{MP}, in terms of a type of $3$-component links called {\em RBG links}. In \cite{MP}, we also put forward five pairs of this type for which $s(K') \neq 0$ (hence $K'$ is not slice), and we could not establish the sliceness of $K$. Soon after, however, Nakamura proved that the knots $K$ in these examples were not slice \cite{Nakamura}.

The idea is still being explored. Dunfield and Gong  \cite{DunfieldGong} found a few more promising examples. The knot $16n68278$ from Figure~\ref{fig:unknown} (c) is of unknown sliceness, but has a $0$-friend that is not slice. The knot $18nh_{00000601}$ from Figure~\ref{fig:unknown} (d) is also of unknown sliceness, and has a $0$-friend that is slice. Interestingly, in this last example the computer program found bands showing that the $0$-friend (which is much larger) is ribbon and hence slice; but could not find this kind of bands for the knot itself. Therefore, either such bands for $18nh_{00000601}$ exist and are hard to find, or the knot gives a counterexample to either SPC4 or the slice-ribbon conjecture.

If disproving SPC4 is a tall order, another (still very ambitious) goal would be to answer Question~\ref{q:definite}  by finding an exotic $\#^n\CP^2$. A strategy is to look for $0$-friends such that one is H-slice in $\#^n\CP^2$ and the other is not. (The Rasmussen invariant still gives an obstruction to H-sliceness in $\#^n\CP^2$, by the results of \cite{MMSW}.) More generally, one can look for $k$-friends (pairs of knots with the same $k$-surgery) such that one is $(-k)$-slice in some  $\#^n\CP^2$ and the other is not; see \cite{Qin} for a discussion. 

\

\bibliographystyle{custom}
\bibliography{biblio2}

\end{document}